\documentclass[a4paper,12pt,twoside]{article}
\usepackage{amsthm,amssymb,amsmath,amscd}
\usepackage{pb-diagram}
\newtheorem{thm}{Theorem}

\newtheorem{lem}{Lemma}
\newtheorem{prop}{Proposition}
\theoremstyle{definition}
\newtheorem{defn}{Definition}
\newtheorem{example}{Example}

\newcommand{\e}{\varepsilon}

\newcommand{\Zset}{\mathbb{Z}}
\newcommand{\Rset}{\mathbb{R}}
\newcommand{\bbc}{\mathbb{C}}
\newcommand{\Cset}{\mathbb{C}}

\newcommand{\bG}{\mathbf{G}}
\newcommand{\btau}{\mathbf{\tau}}

\def\fp{{\mathfrak p}}

\def\cB{{\mathcal B}}
\def\cC{{\mathcal C}}

\def\val{{\rm val}}
\def\Sym{{\rm Sym}}

\def\ie{{\it i.e.,\,}}
\def\cI{{\mathcal I}}

\def\integers{{\mathfrak o}}
\def\ev{{\rm ev}}
\def\odd{{\rm odd}}
\def\RH{{\rm H}}
\def\HP{{\rm HP}}
\def\SL{{\rm SL}}
\def\GL{{\rm GL}}
\def\Aut{{\rm Aut}}
\def\Hom{{\rm Hom}}
\def\End{{\rm End}}
\def\Id{{\rm I}}
\def\Ind{{\rm Ind}}
\def\Res{{\rm Res}}
\def\Mat{{\rm M}}
\def\Nor{{\rm N}}

\def\lcm{{\rm lcm}}

\def\cA{{\mathcal A}}
\def\cR{{\mathcal R}}
\def\cH{{\mathcal H}}
\def\cW{{\mathcal W}}
\def\cL{{\mathcal L}}

\def\fa{{\mathfrak a}}

\def\fR{{\mathfrak R}}
\def\fs{{\mathfrak s}}
\def\ft{{\mathfrak t}}
\def\fA(\fs){{\cA^{\fs}}}
\def\fB{{\mathfrak B}}
\def\fK{{\mathfrak K}}

\def\oA{{\mathfrak A}}
\def\oB{{\mathfrak B}}
\def\oP{{\mathfrak P}}
\def\oQ{{\mathfrak Q}}

\def\toB{{\widetilde{\mathfrak B}}}
\def\um{{\underline{m}\,}}

\def\ZZ{{\mathbb Z}}

\def\fB{{\mathfrak B}}

\def\St{{\rm St}}

\def\Ad{{\rm Ad}}
\def\tJ{{\widetilde J}}
\def\tK{{\widetilde K}}
\def\tM{{\widetilde M}}
\def\tV{{\widetilde V}}

\def\barM{{\bar M}}

\def\barK{{\bar K}}
\def\barN{{\bar N}}
\def\bartau{{\bar \tau}}
\def\ttau{{\widetilde \tau}}
\def\talpha{{\widetilde \alpha}}
\begin{document}
\title{Cycles in the chamber homology of $\GL(3)$}
\author{Anne-Marie Aubert, Samir Hasan and Roger Plymen}
\date{}
\maketitle
\begin{abstract}   Let $F$ be a nonarchimedean
local field and let $\GL(N) = \GL(N,F)$. We prove the existence of
parahoric types for $\GL(N)$. We construct representative cycles
in all the homology classes of the chamber homology of $\GL(3)$.
\end{abstract}

\section{Introduction}

Let $F$ be a nonarchimedean local field and let $G = \GL(N) =
\GL(N,F)$.  The enlarged building $\beta^1 G$ of $G$ is a
polysimplicial complex on which $G$ acts properly. We select a
chamber $C \subset \beta^1 G$. This chamber is a polysimplex, the
product of an $n$-simplex by a $1$-simplex:
\[
C = \Delta_n \times \Delta_1.\]

To this datum we will attach a \emph{homological coefficient
system}, see \cite[p.11]{GM}.  To each simplex $x \in \Delta_n$ we
attach the representation ring $R(G(x))$ of the stabilizer $G(x)$,
and to each inclusion $x \to y$ we attach the induction map:
\[
\Ind_{G(x)}^{G(y)}: R(G(x)) \to R(G(y)).\]

This creates the homology of the simplicial set $\Delta_n$ with
the above coefficient system.   The resulting homology groups are
denoted $h_j(G), 0 \leq j \leq N-1$.

For each point $\fs$ in the Bernstein spectrum $\mathfrak{B}(G)$
(see appendix B) we will select an $\fs$-type $(J,\tau)$.  Here,
$J$ denotes a certain compact open subgroup of $G$, and $\tau$
denotes a certain irreducible smooth representation of $J$, see
\cite{BK3,BK1,BK2}.

The following result is due to Bushnell-Kutzko \cite{BK3, BK1,
BK2}.
\begin{thm}\label{one}  {\rm Existence of types.}
Let $\fs \in \fB(G)$. There exists an $\fs$-type $(J,\tau)$.
\end{thm}

Let $\fs \in \fB(G)$. An $\fs$-type $(J,\lambda)$ will be called
\emph{parahoric} if $J$ is a parahoric subgroup of $G$.

Our first result is the following theorem.

\begin{thm}\label{three} {\rm Existence of parahoric types.}
Let $\fs \in \fB(G)$.  Then there exists a parahoric $\fs$-type
$(J^{\fs},\tau)$.
\end{thm}

The parahoric subgroup $J^{\fs}$ only depends on certain
invariants attached to $\fs$.  For details of these invariants,
see appendix D.

In the proof of Theorem \ref{three}, we have to call upon several
of the technical resources developed by Bushnell-Kutzko.

We now specialize to $\GL(3)$. In this article, we will explicitly
construct representative cycles in \emph{all} the homology classes
in $h_0(G) \oplus h_1(G) \oplus h_2(G)$ when $G = \GL(3)$. This
allows us to compute the chamber homology groups of $\GL(3)$
according to the following formulas:
\[\RH_{\ev}(G;\beta^1 G) = h_0(G) \oplus h_1(G)
\oplus h_2(G) = \RH_{\odd}(G; \beta^1 G).\]

We will demonstrate that each parahoric $\fs$-type
$(J^{\fs},\tau)$ creates finitely many cycles in $h_0(G) \oplus
h_1(G) \oplus h_2(G)$. To prove that all homology classes in
$h_0(G) \oplus h_1(G) \oplus h_2(G)$ are thereby accounted for, we
invoke the $K$-theory of the reduced $C^*$-algebra $\cA: =
C^*_r(G)$. The $K$-theory is torsion-free \cite{P}.

The abelian groups $\RH_{\ev/\odd}(G;\beta^1 G)$ and $K_j(\cA)$
admit compatible Bernstein decompositions, see appendix B.  This
leads, for each $\fs \in \fB(G)$, to the equalities
\begin{eqnarray}
 \textrm{rank} \,\RH_{\ev/\odd}(G; \beta^1 G)^{\fs} = \textrm{rank}
 \,K_0(\cA^{\fs}) = \textrm{rank} \,K_1(\cA^{\fs}).
 \end{eqnarray} The ranks of the finitely generated abelian
groups on the right-hand-side are easily computed (see appendix
C).

\begin{thm}\label{four} Let $G = \GL(3)$, and let $\fs = [M,\sigma]_G$.
  Each parahoric $\fs$-type $(J^{\fs},\tau)$ creates finitely many cycles
in $h_0(G) \oplus h_1(G) \oplus h_2(G)$, and all homology classes
in $h_0(G) \oplus h_1(G) \oplus h_2(G)$ are thereby accounted for.
Quite specifically, we have
\begin{itemize}
\item if $M = \GL(3)$ then \[\RH_{\ev}(G; \beta^1 G)^{\fs} =
\mathbb{Z} = \RH_{\odd}(G; \beta^1 G)^{\fs}\]
 \item if $M = \GL(2)
\times \GL(1)$ then \[\RH_{\ev}(G; \beta^1 G)^{\fs} = \mathbb{Z}^2
= \RH_{\odd}(G; \beta^1 G)^{\fs}\] \item if $M = \GL(1) \times
\GL(1) \times \GL(1)$ then \[\RH_{\ev}(G; \beta^1 G)^{\fs} =
\mathbb{Z}^4 = \RH_{\odd}(G; \beta^1 G)^{\fs}\]
\end{itemize}
\end{thm}

  From this point of view, the types for $\GL(3)$ exceed their
original expectations.  Let $\widehat{\cA^{\fs}}$ denote the dual
of the $C^*$-algebra $\cA^{\fs}$.  This is a compact Hausdorff
space. Since $K$-theory for unital $C^*$-algebras is compatible
with topological $K$-theory of compact Hausdorff spaces, we have
\[ K_j(\fA(\fs)) \cong K^j(\widehat{\fA(\fs)}).\]
Therefore, the $\fs$-type also computes the topological $K$-theory
of the compact space $\widehat{\fA(\fs)}$. The space
$\widehat{\fA(\fs)}$ is precisely the space of all those tempered
representations of $\GL(3)$ which have inertial support $\fs$.

Sections $4-6$ are devoted to a proof of Theorem \ref{three}, and
sections $7-9$ are devoted to a proof of Theorem \ref{four}.

Preliminary work in the direction of Theorem \ref{four} was done
with Paul Baum and Nigel Higson, and recorded in \cite{BHP2}. The
diagrams in \cite{BHP2} are relevant to the present article.  In
\cite{BHP2} all computations were in the \emph{tame} case. We
confront here the general case: this is much more technical. We
require much detailed information in the theory of types; in
particular we need detailed information concerning \emph{compact}
intertwining sets.

We thank the referees for their detailed and constructive
comments.

\section{General results on types}
\label{compactintertwining}

We will collect here some general results on types which will used
in the paper. In this section $G$ denotes the group of $F$-points of an
arbitrary reductive connected algebraic group $\bG$ defined over $F$.

Let $\fR(G)$ denote the category of smooth complex
representations of $G$. Recall that, for each irreducible smooth
representation $\pi$ of $G$, there exists a Levi subgroup $L$ of
a parabolic subgroup $P$ of $G$ and an irreducible
supercuspidal representation $\sigma$ of $L$ such that $\pi$ is
equivalent to a subquotient of the parabolically induced
representation $I_P^G(\sigma)$. The pair $(L,\sigma)$ is unique up
to conjugacy and the inertial class $\fs=[M,\sigma]_G$ (see appendix~B)
is called the \emph{inertial support} of $\pi$.

We have the standard decomposition (see \cite[(2.10)]{Be})
\begin{equation} \label{Bdecomp}
\fR(G) = \prod_{\fs\in\fB(G)}\fR^\fs(G)
\end{equation}
into full sub-categories, where the objects of $\fR^\fs(G)$ are those
smooth representations of $G$ all of whose irreducible subquotients have
inertial support $\fs$.

Let $\fs$ be a point in the Bernstein spectrum of $G$, and let $(J,\tau)$ be an
$\fs$-type, \ie $\tau$ is an irreducible smooth representation of an open
compact subgroup $J$ of $G$ such that for any irreducible smooth representation
$\pi$ of $G$, the restriction of $\pi$ to $J$ contains $\tau$ if and
only if $\pi$ is an object of $\fR^\fs(G)$, \cite[(4.2)]{BK1}.
When $G=\GL(N,F)$, it has been proved \cite{BK3, BK2} that there exists an
$\fs$-type for each point $\fs$ in $\fB(G)$.

\begin{prop} \label{Indtype}
Let $K\supset J$ be an open compact subgroup of $G$ such that
$\alpha:=\Ind_J^K\tau$ is irreducible. Then $(K,\alpha)$ is an
$\fs$-type.
\end{prop}
\begin{proof}
Let $\pi$ be an irreducible smooth representation of $G$.
Using Frobenius reciprocity, we see that
$$\Hom_K\left(\alpha,\Res_K^G(\pi)\right)=
\Hom_J\left(\tau,\Res_J^G(\pi)\right).$$
The result follows from the definition of $\fs$-types.
\end{proof}

\smallskip

Let $J$, $J'$, $K$ be subgroups of $G$ with $J$, $J'$ compact open
and $J\subset K$, $J' \subset K$.  Let $\lambda$, $\lambda'$ be
representations of $J$, $J'$ on finite-dimensional vector spaces
$V$, $V'$.  Let $g\in G$. Then $gJg^{-1} \cap J'$ is a subgroup of
$J'$. We set ${}^g\lambda(x): = \lambda(g^{-1}xg)$.   We define the
$g$-\emph{intertwining vector space} of $(\lambda,\lambda')$ to be
\[\mathcal{I}_g (\lambda,\lambda') = \; \Hom_{gJg^{-1} \cap J'}
({}^g\lambda, \lambda').\] We will write $\cI_g(\lambda) =
\cI_g(\lambda,\lambda)$.

\begin{defn} \label{entrelacement}
\begin{itemize}
\item[(1)]
We say that $g$ {\it intertwines} $\lambda$ if $\cI_g (\lambda) \neq 0$. The
{\it $K$-intertwining set of $\lambda$} is
\[\cI_K(\lambda) = \{g \in K\, :\, \cI_g (\lambda) \neq 0 \}\,\subset\, K.\]
\item[(2)]
We say that $g$ {\it intertwines $\lambda$ and $\lambda'$} if
$\cI_g (\lambda,\lambda')\neq 0$. The {\it $K$-intertwining set of
$\lambda$ and $\lambda'$} is
\[\cI_K(\lambda,\lambda') = \{g \in K\, :\, \cI_g (\lambda,\lambda') \neq 0 \}
\,\subset\, K.\]
\end{itemize}
\end{defn}

In \cite{BK1, BK2, BK3}, the results centre around identification
of the  $G$-intertwining set $I_G(\lambda)$.   In our
applications, we shall need only the $K$-intertwining set
$I_K(\lambda)$ where $K$ is compact.

\smallskip
In order to study the induced representations and their
decomposition into irreducible constituents, we need to use the
Mackey formulas repeatedly.

We assume now that $K$ is open compact. Then
$J$, $J'$ have finite index in $K$.  We have the Mackey
formula:
\begin{equation} \label{Mackey}
\Hom_K \; (\Ind_{J}^K(\lambda), \; \Ind_{J'}^K(\lambda')) \cong \bigoplus
\; \cI_x (\lambda,\lambda')\end{equation} with $x \in J \backslash K /
J'$. If $\lambda = \lambda \cong \lambda'$ then we set $\cI_g
(\lambda) = \cI_g (\lambda,\lambda')$ and we then have the
isomorphism of $\mathbb{C}$-vector spaces
\[{\rm End}_K \; (\Ind_J^K(\lambda)) \cong \bigoplus \; \cI_x (\lambda)\]
with $x \in J \backslash K / J$.


The following is an immediate consequence: we will use this result
repeatedly.

\begin{prop} \label{IK} If $\cI_K(\lambda) = J$ then $\Ind_J^K
(\lambda)$ is irreducible.
\end{prop}

We will use the following immediate result.

\begin{prop} \label{IKequiv} If $\Ind_{J}^K(\lambda)$ and
$\Ind_{J'}^K(\lambda')$ are irreducible, then
\[\Ind_{J}^K(\lambda)\cong\Ind_{J'}^K(\lambda') \Longleftrightarrow
\cI_K(\lambda,\lambda') = JyJ'\]\ for some element
$y$.\end{prop}

\begin{prop} \label{fromtype}
Let $(J^\fs,\tau^\fs)$ be an $\fs$-type, $(J^{\fs'},\tau^{\fs'})$ be a
$\fs'$-type with $\fs$, $\fs'$ in $\fB(G)$, $\fs \neq \fs'$. Let $J$ be a
compact open subgroup of $G$ such that $J^\fs \subset J$, $J^{\fs'}
\subset J$. Then we have
\[\dim_{\Cset}\Hom_J(\Ind_{J^\fs}^{J}\,\tau^\fs,\Ind_{J^{\fs'}}^{J}\,\tau^{\fs'})
= 0.\]
\end{prop}
\begin{proof} From the Mackey formula~(\ref{Mackey}), it is
equivalent to prove that $\cI_J(\tau^\fs,\tau^{\fs'})=0$. The
proof of the equivalence of (i) and (ii) of
\cite[Theorem~9.3.a]{BK1} shows that $\fs=\fs'$ if and only if
$\cI_G(\tau^\fs,\tau^{\fs'})\ne 0$. The result follows.
\end{proof}

\smallskip

Let $J$ be a compact open subgroup of $G$ and
$(\tau,\cW)$ be an irreducible smooth representation of $J$. Let
$(\tau^\vee,\cW^\vee)$ be the contragredient representation of
$(\tau,\cW)$.

For any subgroup $K$ of $G$, let $\cH(K,\tau)$ denote the space of
compactly supported functions $f\colon K\to\End_\Cset(\cW^\vee)$ such that
$f(j_1kj_2)=\tau^\vee(j_1)f(k)\tau^\vee(j_2)$, for any $j_i\in J$, $k\in
K$. The standard convolution operation gives $\cH(K,\tau)$ the structure
of an associative unital $\Cset$-algebra.

\smallskip

Let $M$ be a Levi subgroup of $G$, and let $(J_M,\tau_M)$ be a
$\ft$-type, with $\ft:=[M,\sigma]_M$ a (supercuspidal) point of the
Bernstein spectrum of $M$.

\smallskip

We recall from
\cite[Definition~8.1]{BK1} that the pair $(J,\tau)$ is a \emph{$G$-cover
of $(J_M,\tau_M)$} if $J\cap M=J_M$ and $\tau_{|J_M}\cong \tau_M$, and
if the following conditions hold for every parabolic
subgroup $P$ of $G$ with Levi subgroup $M$:
\begin{itemize}
\item[(1)]
$(J,\tau)$ it is \emph{decomposed with respect to $(M,P)$}, that is, $J$
admits the Iwahori decomposition: $$J=J\cap U\cdot J_M\cdot J\cap \overline U,$$
and the groups $J\cap U$, $J\cap \overline U$ are both contained in the kernel
of $\tau$ (here $U$, $\overline U$ denote the unipotent radicals of $P$
and of its opposite parabolic subgroup, respectively),
\item[(2)]
there exists an invertible element of $\cH(G,\tau)$ supported on a double
coset $Jz_PJ$, where $z_P$ is a central element in $M$, which is strongly
\emph{$(P,J)$-positive} in the sense of \cite[Definition~(6.16)]{BK1}.
\end{itemize}

\smallskip

The group $\Psi(M)$ of unramified quasicharacters of $M$ has the structure
of a complex torus.
The action (by conjugation) of $\Nor_G(M)$ on $M$ induces an action of
$W(M):=\Nor_G(M)/M$ on $\fB(M)$. Let $W_\ft$ denote the stabilizer of
$\ft=[M,\sigma]_M$ in $W(M)$. Thus $W_\ft=N_\ft/M$, where
\begin{equation} \label{Nor_t}
\Nor_\ft=\left\{n\in\Nor_G(M)\,:\,{}^n\sigma\cong\nu\sigma, \;\text{ for
some $\nu\in \Psi(M)$}\right\}\end{equation}
denotes the $\Nor_G(M)$-normalizer of $\ft$.

\smallskip
We will need the following Proposition which gives a bound for the
compact intertwining.

\begin{prop} {\rm \cite{BK4}} \label{Intbound}
We assume here that $G=\GL(N,F)$.
Let $M$ be a Levi subgroup of $G$, let $(J,\tau)$ be a $G$-cover
of a $\ft$-type, with $\ft=[M,\sigma]_M$ a point of the
Bernstein spectrum of $M$, and let $K$ be a compact subgroup of
$G$ which contains $J$. Let $t$ denote the number of double
classes $J\backslash K/J$ which intertwine $\tau$. Then
$$t\,\le\,|W_{\ft}|.$$
\end{prop}
\begin{proof} It is a classical result that $t$ is bounded by the dimension of
$\mathcal{H}(K,\tau)$.
The hypotheses on the the supercuspidal representation $\sigma$ which are
listed in \cite[\S 1.3]{BK4} are identical to those listed in \cite[(5.5)]{BK1}.
Since $G=\GL(N,F)$, it follows from \cite[Comments~(b) and (d) on (5.5)]{BK1}
that these hypotheses are satisfied, and so we can
apply \cite[Theorem 1.5(ii)]{BK4}.   We infer that
\[\dim_{\Cset} \cH(K,\tau) \leq |W_{\ft}|.\]
\end{proof}
\section{Chamber homology groups}

Let $\integers_F$ denote the ring of integers of $F$, let
$\varpi=\varpi_F$ be a uniformizer in $F$, and
$\fp_F=\varpi_F\integers_F$ denote the maximal ideal of $\integers_F$.
We set
\[\Pi=\Pi_N=\left(\begin{matrix}
0&\Id_{N-1}\cr \varpi&0\end{matrix}\right).\] Let $s_0$, $s_1$,
$\ldots$, $s_{N-1}$ denote the standard involutions in $G$: $s_i$
denote the matrix in $G$ of the transposition $i\leftrightarrow
i+1$, that is,
\[s_i=
\left(\begin{matrix}
\Id_{i-1}&&&\cr
&0&1&\cr
&1&0&\cr
&&&\Id_{N-i-1}
\end{matrix}\right),\]
for every $i\in\{1,\ldots,N-1\}$, and $s_0=\Pi s_1\Pi^{-1}$.

The finite Weyl group is $W_0=<s_1,s_2,\ldots,s_{N-1}>$, and the
affine Weyl group is given by $W=<s_0,s_1,\ldots,s_{N-1}>$.

We set
\[
\cR(g) = \Pi^{-1}g\Pi
\]
with $g \in G$, so that $\cR^N = 1$.

\smallskip

We will use repeatedly, and without further comment, the fact that
induction commutes with conjugation: in particular conjugation by
$\Ad \, \Pi^i$, $1\le i\le N-1$.  We will use this in the following form:
\begin{equation} \label{RI}
{\cR}^{-1}(\Ind_{\cR H}^{\cR G}(\cR \alpha)) \cong
\Ind_H^G(\alpha).
\end{equation}

\smallskip
Note that
\[
\cR(s_i) = s_{i+1},\;\;\text{
with $i = 0$, $1$, $\ldots$, $N-1$  $\mod N$.}\]

The extended affine Weyl group is given by $\widetilde{W} = W
\rtimes <\Pi> $. We observe that
\begin{equation} \label{WK}
\widetilde{W}\cap\GL(N,\integers_F)=W_0.
\end{equation}The standard Iwahori subgroup is
$$ I = \left( \begin{array}{cccc}
\integers_F^{\times} &\integers_F &\cdots &\integers_F\\
\;\fp_F & \ddots &\ddots &\vdots\\
\vdots &\ddots&\ddots &\integers_F \\
\;\fp_F &\cdots&\;\fp_F& \integers_F^{\times}
\end{array} \right).
$$
Let $A$ be the apartment attached to the diagonal torus and let
$\Delta$ denote the unique chamber of of $A$ which is stabilized
by $<\Pi>I$. We index the vertices $L_0$, $L_1$, $\ldots$,
$L_{N-1}$ of $\Delta$ in such a way that
\begin{itemize}
\item[$\bullet$] $s_i\Delta$ is the unique chamber of $A$ which is
adjacent to $\Delta$ and such that $s_i\Delta\cap\Delta$ is the
$(N-2)$-simplex $\{L_0,\ldots,L_{N-1}\}\backslash\{L_i\}$;
\item[$\bullet$] $\cR(L_i)=L_{i+1}$ with $i = 0$, $1$, $\ldots$,
$N-1$ $\mod N$.
\end{itemize}
The $L_i$ are the maximal standard parahoric subgroups of $G$,
\[L_i=I<s_0,s_1, \ldots,s_{i-1},s_{i+1},\ldots,s_{N-1}>I\,=\,
\cR^i(L_0),\]
and $L_0=\GL(N,\integers_F)$.

The stabilizers of the facets of dimension $N-1$ of $\Delta$ are
$J_0$, $J_1$, $\ldots$, $J_{N-1}$, where
\[J_i=I<s_i>I.\]
Each parahoric subgroup of $G$ is defined by a facet of the building and the
standard parahoric subgroups are the
\[J_S=I<s_j\,:\,j\in S>I,\]
where $S$ is any subset of $\{0,1,\ldots,N-1\}$ mod. $N$, \cite[p.  118]{Ron}.

Hence,
$I=J_{\emptyset}$, $J_i=J_{\{i\}}$, $L_i=J_{\{0,1,\ldots,i-1,i+1,\ldots,N-1\}}$.

The enlarged building $\beta^1 G$ is labellable, that is, there
exists a simplicial map $\ell\colon\beta^1 G\to\Delta$, which
preserves the dimensions of the simplices. The labelling is
unique, up to the automorphisms of $\Delta$. It allows us to fix
an orientation of the simplices: one defines an incidence number
$<\eta:\sigma>$ between an arbitrary facet
$\eta=(\eta_0,\ldots,\eta_{i-1})$ of dimension $i$ and any facet
$\sigma=(\sigma_0,\ldots,\sigma_i)$ of dimension $i+1$ which
contains $\eta$, as follows
\[<\eta:\sigma>=(-1)^i\;\;\text{if
$\{\ell(\eta_0),\ldots,\ell(\eta_{i-1})\}\backslash
\{\ell(\sigma_0),\ldots,\ell(\sigma_i)\}=i$.}\]
\smallskip

The chamber homology groups are obtained by totalizing the
bicomplex $C_{**}$
$$
\begin{array}{lllllllllllll}
0 & \longleftarrow & C_0  &\longleftarrow &\cdots&
\longleftarrow &C_i&\longleftarrow &\cdots&\longleftarrow & C_{N-2}
 &\longleftarrow & C_{N-1}\\
  &                & \downarrow &                && &
\downarrow &                &
& &\downarrow  &              &
 \downarrow\\
0 & \longleftarrow & C_0  & \longleftarrow  &\cdots &
\longleftarrow &C_i&\longleftarrow &\cdots&\longleftarrow
& C_{N-2}& \longleftarrow & C_{N-1}\\
\end{array}
$$
in which the chains are as follows:
\begin{equation} \label{Ci}
C_i=\bigoplus_{S\subset\{0,1,\ldots,N-1\}\atop |S|=N-1-i}R(J_S)
\end{equation}
and each vertical map is given by $1 - \Ad\,\Pi$.
In particular, we have
\begin{itemize}
\item $C_0 = R(L_0) \oplus R(L_1) \oplus \cdots \oplus R(L_{N-1})$,
\item $C_{N-2} = R(J_0) \oplus R(J_1) \oplus \cdots \oplus R(J_{N-1})$,
\item $C_{N-1} = R(I)$.
\end{itemize}

We will write an arbitrary element $v$ in $C_i$ as a
$\left(\begin{smallmatrix}N\cr i\end{smallmatrix}\right)$-uple
$[\eta]$.
Once an orientation has been chosen, the differentials are as
follows: if $v \in C_i$ then
\[ \partial(v)\,= \,\sum_{\eta\subset\sigma\atop{\dim\eta=i}}
(-1)^{\langle\eta:\sigma\rangle}\Ind^{G(\eta)}_{G(\sigma)}[\eta]\,\in
C_{i-1}.
\]
In particular:
\begin{itemize}
\item if $v \in C_{N-1}$ then $\partial(v) =
(\Ind_I^{J_0}(v),\Ind_I^{J_1}(v),\ldots,\Ind_I^{J_{N-1}}(v))$,
\item if $v \in C_0$ then $\partial(v) = 0$.
\end{itemize}

When $G=\GL(3)$, if $v =(v_0,v_1,v_2) \in C_1$ then $\partial(v)$ equals
\[(\Ind_{J_2}^{L_0}(v_2)
- \Ind_{J_1}^{L_0}(v_1),\Ind_{J_0}^{L_1}(v_0) -
\Ind_{J_2}^{L_1}(v_2),- \Ind_{J_0}^{L_2}(v_0) +
\Ind_{J_1}^{L_2}(v_1)),\]
and, in the chain complex
\[
0 \longleftarrow C_0 \longleftarrow C_1 \longleftarrow
\longleftarrow C_2\longleftarrow 0,
\]
we have that $v$ is a $1$-cycle if and only if
\[
\Ind_{J_2}^{L_0}(v_2) = \Ind_{J_1}^{L_0}(v_1),
\Ind_{J_0}^{L_1}(v_0) = \Ind_{J_2}^{L_1}(v_2),
\Ind_{J_0}^{L_2}(v_0) = \Ind_{J_1}^{L_2}(v_1),\] \ie if and only
if the $1$-chain $(v_0,v_1,v_2)$ is \emph{vertex compatible}. Note
that a true representation in $R(I)$ can never be a $2$-cycle; on
the other hand, each $0$-chain is a $0$-cycle.


When we totalize the bicomplex we obtain the chain complex
\[
0 \longleftarrow C_0 \longleftarrow C_0 \oplus C_1 \longleftarrow\cdots
\longleftarrow
C_{i-1} \oplus  C_i \longleftarrow C_i\oplus C_{i+1}
\longleftarrow\cdots\longleftarrow C_{N-1} \longleftarrow 0
\]

\begin{defn} The homology groups of this totalized complex are the
chamber homology groups, as in \cite{BHP2}.
\end{defn}

To each point $\fs \in \mathfrak{B}(G)$ we will associate a
sub-bicomplex $C_{**}(\fs)$:

$$
\begin{array}{lllllllllllll}
0 & \longleftarrow & C_0(\fs)  &\longleftarrow &\cdots&
\longleftarrow &C_i(\fs)&\longleftarrow &\cdots&\longleftarrow &  C_{N-1}(\fs)\\
  &                & \downarrow &                && &
\downarrow &                & & &\downarrow  &              &
 \\
0 & \longleftarrow & C_0(\fs)  & \longleftarrow  &\cdots &
\longleftarrow &C_i(\fs)&\longleftarrow &\cdots& \longleftarrow & C_{N-1}(\fs)\\
\end{array}
$$
in which each vertical differential is  $0$.  By an
\emph{invariant chain} we shall mean a chain invariant under $\Ad
\, \Pi$.

Let $\fs$ be a point in $\mathfrak{B}(G)$ with $\fs = [M,\sigma]_G$.
We recall that $W(M)$ denotes the group $\Nor_G(M)/M$.
We take for $M$ a standard Levi subgroup of
$G$, isomorphic to $\GL(N_1)\times\cdots\times\GL(N_r)$, with
$(N_1\ge N_2\ge\cdots\ge N_r)$ a partition of $N$.

Given a point $\fs \in \mathfrak{B}(G)$, fix an $\fs$-type
$(J,\tau)$. Such types exist \cite{BK3, BK1, BK2}.  There exists a
parahoric subgroup $J^\fs$ containing $J$ such that
$(J^\fs,\alpha:=\Ind_J^{J^\fs}\tau)$ is also an $\fs$-type
(see Theorems~\ref{simplet}, \ref{homogenet}, \ref{ssinclusion}).

Then
\begin{itemize}
\item induce (if possible) each element in the orbit $W(M)\cdot
\alpha$ to the standard parahoric subgroups containing $J^{\fs}$,
and rotate, \ie apply $\cR$, $\ldots$, $\cR^{N-1}$, \item  take
the free abelian groups generated by all the irreducible
components which arise in this way.
\end{itemize}

Each of our sub-complexes $C_{**}(\fs)$ will come from some or all
of this data. All the chain groups in $C_{**}(\fs)$ are finitely
generated free abelian groups and comprise invariant chains. The
homology groups of the chain complex
\[
0 \longleftarrow C_0(\fs) \longleftarrow C_1(\fs) \longleftarrow
\cdots\longleftarrow C_{N-1}(\fs) \longleftarrow 0
\]
will be denoted $h_*(\fs)$.  We call this the \emph{little
complex}.

When we totalize the associated bicomplex $C_{**}(\fs)$ we obtain
the chain complex
\[
0 \longleftarrow C_0(\fs)
\longleftarrow\cdots\longleftarrow
C_{i-1}(\fs) \oplus  C_i(\fs) \longleftarrow C_i(\fs)\oplus C_{i+1}(\fs)
\longleftarrow\cdots\longleftarrow C_{N-1}(\fs) \longleftarrow 0
\]

The following lemma will speed up our calculations.
\begin{lem} \label{Homology Lemma}
The homology groups $H_*(\fs)$ of this complex are given by
\[
H_0(\fs) = h_0(\fs),\;\;\; H_N(\fs) = h_{N-1}(\fs)
\]
\[H_{i+1}(\fs) = h_i(\fs) \oplus h_{i+1}(\fs),  \;\; 0 \leq i \leq
N-2
\]
\[
H_{\ev}(\fs) = h_0(\fs) \oplus h_1(\fs) \oplus \cdots \oplus
h_{N-1}(\fs) = H_{\odd}(\fs)
\]
The even (resp. odd) chamber homology is precisely the total
homology of the little complex.
\end{lem}

\begin{proof}  This is a direct consequence of the fact that each vertical differential
in the bicomplex $C_{**}(\fs)$ is $0$.
\end{proof}


\section{Lattice chains and lattice sequences}

Let $V$ be an $F$-vector space of dimension $N$.
We recall from \cite[Def.~2.1]{BK2} that a \emph{lattice sequence} is
a function $\Lambda$ from $\Zset$ to the set of $\integers_F$-lattices in
$V$ such that
\begin{itemize}
\item
$i\ge j$ implies $\Lambda(i)\le\Lambda(j)$;
\item
there exists $e=e(\Lambda)\in\Zset$, $e\ge 1$, such that
$\Lambda(i+e)=\fp_F\,\Lambda(i)$ for any $i\in\Zset$.
\end{itemize}
The integer $e$ is uniquely determined, and is called the
\emph{period} of $\Lambda$. We have $e\le N$.

A lattice sequence which is injective as a function is called
\emph{strict}.
We will put
\begin{equation} \label{aLambda}
\fa_n(\Lambda):=\left\{a\in
A\,:\,a\Lambda(m)\subset\Lambda(m+n),\;m\in\Zset\right\},\;\;n\in\Zset.
\end{equation}

The concept of lattice sequence generalizes the notion of lattice chain:
as defined in \cite[(1.11)]{BK3}, a \emph{lattice chain} in $V$ is a
set $\cL=\left\{L_i\,:\,i\in\Zset\right\}$ of $\integers_F$-lattices $L_i$
in $V$ such that
\begin{itemize}
\item
$L_i\supset L_{i+1}$, $L_i\ne L_{i+1}$, for any $i\in\Zset$;
\item there exists $e=e(\cL)\in\Zset$ such that $L_{i+e}=\fp_F\,L_i$, for any
$i\in\Zset$.
\end{itemize}
The integer $e$ is uniquely determined, and is called the
\emph{period}
of $\cL$.

Let $k_F$ denote the residue field of $F$. For each $i$, the quotient
$L_i/L_{i+1}$ is a $k_F$-vector space. Write
\begin{equation} \label{di}
d_i=d_i(\cL):=\dim_{k_F}(L_i/L_{i+1}).
\end{equation}
The function $d(\cL)\colon i\mapsto d_i$, $i\in\Zset$, is periodic of
period dividing $e$, and we have
\begin{equation} \label{sum_di}
\sum_{i=0}^{e-1}d_i=N.
\end{equation}

To each lattice chain $\cL$ is attached a strict lattice sequence
$\Lambda_\cL$ defined by $\Lambda_\cL(i):=L_i$, for $i\in\Zset$.
In the opposite direction, to each lattice sequence $\Lambda$ is attached
a lattice chain $\cL_\Lambda$  defined by
\begin{equation} \label{LL}
\cL_\Lambda:=\left\{\Lambda(i)\,:\,i\in\Zset\right\}.
\end{equation}

\smallskip

As in \cite[\S 2.6]{BK2}, we extend a lattice sequence $\Lambda$ to a function
on the real line $\Rset$ by setting
\begin{equation} \label{reel}
\Lambda(x):=\Lambda(\lceil x\rceil),
\;\;\text{ $x\in\Rset$,}
\end{equation} where
$\lceil x\rceil$ is the integer defined
by the relation $\lceil x\rceil-1<x\le \lceil x\rceil$.

\smallskip

Let $\Lambda$ be a lattice sequence in $V$ and let $m$ be a positive integer.
Then the function $m\Lambda$ from $\Zset$ to the set of $\integers_F$-lattices
in $V$ defined by
$$(m\Lambda)(i):=\Lambda(i/m),\;\;\text{ for any $i\in\Zset$,}$$
is a lattice sequence in $V$ with period $m\,e(\Lambda)$, and we have
\begin{equation} \label{mL}
(m\Lambda)(i)=\begin{cases}\Lambda(i/m)&\text{if $m$ divides $i$,}\cr
\Lambda(1+[i/m])&\text{otherwise,}
\end{cases}
\end{equation}
and
$(m\Lambda)(x)=\Lambda(x/m)$, for all $x\in\Rset$ (see
\cite[Prop.~2.7]{BK2}).

\smallskip

If we have a lattice sequence $\Lambda$ in $V$ and an integer $t$, we can
define a lattice chain $\Lambda+t$ by
\begin{equation} \label{Lt}
(\Lambda+t)(i):=\Lambda(i+t), \;\;\text{for
any $i\in\Zset$.}
\end{equation}

\smallskip

Let $m$ be a positive integer, and let $V^1$, $V^2$, $\ldots$, $V^m$ be $m$
finite-di\-men\-sio\-nal $F$-vector spaces.
Let $\Lambda^1$, $\Lambda^2$, $\ldots$, $\Lambda^m$ be $m$ lattices sequences
in $V$,
with periods $e_1$, $e_2$, $\ldots$, $e_m$, respectively.
We denote by $\Lambda=\Lambda^1\oplus\cdots\oplus\Lambda^m$ the
\emph{direct sum} of
$\Lambda^1$, $\ldots$, $\Lambda^m$: we recall from \cite[\S 2.8]{BK2} that
$\Lambda$ is defined by
\begin{equation} \label{seq_add}
\Lambda(ex)=\Lambda^1(e_1x)\oplus\cdots\oplus\Lambda^m(e_mx),\;\;\text{
for each $x\in\Rset$,
where $e=\lcm\{e_1,\ldots,e_m\}$.}\end{equation}

\smallskip

The following example occurs in the construction of \cite[\S 7.2]{BK2}.
See also \cite[Example 2.8]{BK2}.

\begin{example} \label{Exemple}
We assume given $m$ lattice chains $\cL^1$, $\cL^2$, $\ldots$, $\cL^m$ in $V^1$,
$V^2$, $\ldots$, $V^m$, respectively, of same period $e$.
We define a lattice chain
$$\cL=\left\{L_i\,:\,i\in\Zset\right\}$$
in $V$ of period $me$ by setting
$$L_{mj+k}:=L_j^1\oplus L_j^2\oplus\cdots\oplus L_j^{m-k}\oplus
L_{j+1}^{m-k+1}\oplus\cdots\oplus L_{j+1}^m,$$
any $j\in\Zset$ and $0\le k\le m-1$.
Using (\ref{mL}), (\ref{Lt}), we obtain
$$\Lambda_\cL=(m\Lambda^1-m+1)\oplus\cdots\oplus
(m\Lambda^{m-k}-k)\oplus\cdots\oplus (m\Lambda^{m-1}-1)\oplus
m\Lambda^m.$$
\end{example}

\subsection{Addition of lattice chains} \label{addition_procedure}


Let $A:=\End_F(V)$ and let $E/F$ be a subfield of $A$. We denote
by $\integers_E$ the discrete valuation ring in $E$, by $k_E$ its residue
field, and by
$e(E|F)$ the ramification degree of $E/F$.

Let $V^1$, $V^2$, $\ldots$, $V^m$ be $m$
finite-di\-men\-sio\-nal $F$-vector spaces of dimensions $N_1$, $N_2$,
$\ldots$, $N_m$, respectively. We assume that the field $E$ preserves the
spaces $V^i$. We may consider each $V^l$ as a $E$-vector
space of dimension $N_l/[E:F]$.

Let $\cL^1$, $\cL^2$, $\ldots$, $\cL^m$ be $m$ $\integers_E$-lattice chains
in the $E$-vector spaces $V^1$, $V^2$, $\ldots$, $V^m$, respectively,
of period $e_1'$, $e_2'$, $\ldots$, $e_m'$, respectively.

\subsubsection{First addition procedure} \label{first}

We define first an $\integers_E$-lattice chain
$\cL^1+\cL^2=\{L_j^{[1,2]}\,:\,j\in\Zset\}$ in $V^1\oplus V^2$ of period
$e_1'+e_2'$ by
\begin{equation} \label{L12}
L_i^{[1,2]}:=\begin{cases}L_0^1\oplus L_i^2,&\text{if $0\le i\le e_2'-1$}\cr
L_{i-e_2'}^1\oplus L_{e_2'}^2,&\text{if $e_2'\le i\le e_1'+e_2'-1$.}
\end{cases}\end{equation}
Then let $\cL^1+\cL^2+\cL^3=\{L_i^{[1,3]}\,:\,i\in\Zset\}$
be the $\integers_E$-lattice chain
$(\cL^1+\cL^2)+\cL^3$ (which is the same as $\cL^1+(\cL^2+\cL^3)$).
By applying~(\ref{L12}) to the two $\integers_E$-lattice chains
$\cL^1+\cL^2$ and $\cL^3$, we get
\[L_i^{[1,3]}:=\begin{cases}L_0^{[1,2]}\oplus L_i^3,&\text{if $0\le i\le
e_3'-1$}\cr
L_{i-e_3'}^{[1,2]}\oplus L_{e_3'}^2,&\text{if $e_3'\le i\le (e_1'+e_2')+e_3'-1$,}
\end{cases}\]
that is, by using~(\ref{L12}),
\begin{equation} \label{L123}
L_i^{[1,3]}:=\begin{cases}L_0^1\oplus L_2^0\oplus L_i^3,&\text{if $0\le i\le
e_3'-1$}\cr
L_0^1\oplus L_{i-e_3'}^2\oplus L_{e_3'}^2,&\text{if $e_3'\le i\le e_2'+e_3'-1$}\cr
L_{i-e_3'-e_2'}^1\oplus L_{e_2'}^2\oplus L_{e_3'}^2,&\text{if $e_2'+e_3'\le i\le
e_1'+e_2'+e_3'-1$.}
\end{cases}\end{equation}
Using this procedure, we finally obtain an $\integers_E$-lattice
chain
\begin{equation} \label{addedchains}
\cL^1+\cdots+\cL^m=\cL^{[1,m]}:=
\left\{L^{[1,m]}_i\,:\,i\in\Zset\right\}
\end{equation} of period
$e_1'+e_2'+\cdots+e_m'$.
We have
\begin{equation} \label{addedchainsII}
L^{[1,m]}_i=L_{0}^1\oplus\cdots\oplus L_{0}^{j-1}\oplus
L_{k}^j\oplus L_{e_{j+1}'}^{j+1}\oplus\cdots\oplus L_{e_m'}^m,
\end{equation}
for $i=e_{j+1}'+\cdots+e_m'+k$ with $1\le j\le m$ and $0\le k\le e_j'-1$.

We will need the $\integers_E$-lattice chain $\cL^1+\cdots+\cL^m$ in the
special case when
$e_1'=\cdots=e_m'=1$. In that case, the equation~(\ref{addedchainsII}) becomes
\begin{equation} \label{sumchain}
L^{[1,m]}_i=L_{0}^1\oplus\cdots\oplus
L_0^{m-i-1}\oplus L_0^{m-i}\oplus L_1^{m-i+1}\oplus\cdots\oplus L_1^m,
\end{equation}
for each $i\in\{0,1,\ldots,m-1\}$.

Since the $\integers_E$-lattice chains $\cL^1$, $\ldots$, $\cL^m$ all have
period $1$, we have $\fp_E^j L^l_0=L_{j}^l$, for each $l\in\{1,\ldots,m\}$ and
each $j\in\Zset$. Hence, since the $\integers_E$-lattice chain
$\cL^{[1,m]}$ is of period $m$, we have
\begin{equation} \label{chaineBK}
L^{[1,m]}_{mj+k}=\fp_E^jL_k^{[1,m]}=
L_{j}^1\oplus\cdots\oplus L_{j}^{m-k-1}\oplus
L_{j}^{m-k}\oplus L_{j+1}^{m-k+1}\oplus\cdots\oplus L_{j+1}^m,
\end{equation}
for each $j\in\Zset$ and each $k\in\{0,1,\ldots,m-1\}$.

Then we have
$$L^{[1,m]}_{mj+k+1}= L_{j}^1\oplus\cdots\oplus L_{j}^{m-k-1}
\oplus L_{j+1}^{m-k}\oplus L_{j+1}^{m-k+1}\oplus\cdots\oplus L_{j+1}^m.
$$
It follows that
$$L^{[1,m]}_{mj+k}/L^{[1,m]}_{mj+k+1}\cong L_j^{m-k}/L_{j+1}^{m-k},$$
for each $j\in\Zset$ and each $k\in\{0,1,\ldots,m-1\}$.

Hence, setting $d^l:=d(\cL^l)$ for any $1\le l\le m$, we obtain
\begin{equation} \label{dL1m}
d_{mj+k}(\cL^{[1,m]})=d^{m-k}_j=d^{m-k}_0=\dim_E(V^{m-k})=N_{m-k}/[E:F],
\end{equation}
by (\ref{sum_di}), since $e_{m-k}=1$.

\smallskip

We may consider each $\cL^l$, $1\le l\le m$, as an $\integers_F$-lattice chain
in the $F$-vector space $V$, of period $e(E|F)$
(see \cite[(1.2.4)]{BK3}). Then $\cL^1+\cdots+\cL^m$, viewed as an
$\integers_F$-lattice chain, has
period $m\,e(E|F)$ (by \cite[(1.2.4)]{BK3}) and the equation~(\ref{chaineBK})
shows that it is the same as the chain $\cL$ considered in the
Example~\ref{Exemple}.

\subsubsection{Second addition procedure} \label{second}

We keep assuming $e_1'=\cdots=e_m'=1$, and we will now consider the
$\integers_E$-lattice chain
\[\cL^m+\cdots+\cL^1=
\left\{L_i^{[m,1]}\,:\,i\in\Zset\right\}
.\]
We have
\begin{equation} \label{goodchain}
L^{[m,1]}_{mj+k}=
L_{j+1}^1\oplus\cdots\oplus L_{j+1}^{k}\oplus
L_{j}^{k+1}\oplus\cdots\oplus L_{j}^m,
\end{equation}
for each $j\in\Zset$ and each $k\in\{0,1,\ldots,m-1\}$.

It gives
\begin{equation} \label{goodquotient}
L^{[m,1]}_{mj+k}/L^{[m,1]}_{mj+k+1}\cong L_j^{k+1}/L_{j+1}^{k+1},
\end{equation}
for each $j\in\Zset$ and each $k\in\{0,1,\ldots,m-1\}$.
Hence we obtain
\begin{equation} \label{dLm1}
d_{mj+k}(\cL^{[m,1]})=d^{k+1}_j=d^{k+1}_0=\dim_E(V^{k+1})=N_{k+1}/[E:F],
\end{equation}
which in particular does not depend on $m$, in contrast with
$d_{mj+k}(\cL^{[1,m]})$.

\smallskip

As before, we may consider each $\cL^l$, $1\le l\le m$, as an
$\integers_F$-lattice chain in the $F$-vector space $V$, of period $e(E|F)$.
Then $\cL^m+\cdots+\cL^1$, viewed as an
$\integers_F$-lattice chain, has period $m\,e(E|F)$.

\section{Hereditary $\integers_F$-orders} \label{orders}

To any $\integers_F$-lattice chain $\cL=\{L_i\}$ in $V$ is attached the
following sequence of $\integers_F$-lattices in $A$
\[\End_{\integers_F}^n(\cL):=
\left\{x\in A\,:\,x L_i\subset L_{i+n}, \;i\in\Zset\right\},\]
for each $n\in\Zset$. In particular, $\oA=\oA(\cL):=\End_{\integers_F}^0(\cL)$ is
an hereditary $\integers_F$-order in $A$, and $\oP:=\End_{\integers_F}^1(\cL)$ is
the Jacobson radical of $\oA$. We will set
\begin{equation} \label{UA}
U(\oA):=\oA^\times\quad\text{and}\quad U^n(\oA):=1+\oP^n,\;\; \text{for
$n\ge 1$.}
\end{equation}
We put
\begin{equation} \label{normalizer}
\fK(\oA): =\left\{g\in \Aut_F(V)\,:\,g^{-1}\oA g=\oA\right\}.
\end{equation}

\smallskip

\begin{defn} \label{Ordre}
For any partition $(N_1,N_2,\ldots,N_r)$ of $N$, we denote by
\[\oA(N_1,N_2,\ldots,N_r)\] the subset of $\Mat_{N}(F)$ consisting of the
matrices of the following form: the $(i,j)$-block has dimension $N_i\times
N_j$, $1\le i,j\le r$, and its entries lie in $\integers_F$ if $i\le j$,
in $\fp_F$ otherwise.
Pictorially,
\[
\oA(N_1,N_2,\ldots,N_r)= \left( \begin{array}{cccc}
\integers_F &\integers_F &\cdots &\integers_F\\
\;\fp_F & \ddots &\ddots &\vdots\\
\vdots &\ddots&\ddots &\integers_F \\
\;\fp_F &\cdots&\;\fp_F& \integers_F
\end{array} \right).
\]
\end{defn}
\smallskip

Let $e:=e(\cL)$ and $d_i:=d_i(\cL)$.
For each $i\in\{0,1,\ldots,e-1\}$, we choose elements $v_{i,h}\in L_i$,
$1\le h\le d_i$ such that the cosets $v_{i,h}+L_{i+1}$ form a basis of the
$k_F$-space $L_i/L_{i+1}$. Then
\[(v_{e-1,1}, v_{e-1,2}, \ldots, v_{e-1,d_{e-1}},
v_{e-2,1}, v_{e-2,2}, \ldots, v_{e-2,d_{e-2}},
\ldots, v_{0,1}, v_{0,2}, \ldots, v_{0,d_{0}}) \]
is an $F$-basis of $V$.
If we use this basis to identify $A$ with the matrix algebra
$\Mat_{N}(F)$, then $\oA$ becomes identified with
$\oA(d_0,d_2,\ldots,d_{e-1})$.

\smallskip

Now let $V^1$, $V^2$, $\ldots$, $V^m$ be $m$ finite-di\-men\-sio\-nal $F$-vector
spaces as in sections~\ref{first}, \ref{second}, and let $\cL^1$, $\cL^2$,
$\ldots$, $\cL^m$ be $m$ $\integers_E$-lattice chains in
$V^1$, $V^2$, $\ldots$, $V^m$, all of period $1$.
We put $$\oA^{[m,1]}:=\oA(\cL^{[m,1]}),$$
where $\cL^{[m,1]}$ is defined as in~(\ref{goodchain}).

Let $(m_1,\ldots,m_r)$ be a partition of $m$. For each
$i\in\{1,\ldots,r\}$, we set $\um_{i-1}:=m_1+\cdots+m_{i-1}$,
\[\cL^{[\um_i,\um_{i-1}+1]}:=\cL^{\um_i}+\cL^{\um_i-1}+\cdots+
\cL^{\um_{i-1}+2}+\cL^{\um_{i-1}+1},\]
and
\[\oA^{[m_i,m_{i-1}+1]}:=\oA(\cL^{[\um_i,\um_{i-1}+1]}).\]
We set $m_0:=0$.
For each $i\in\{1,\ldots,r\}$, we define $V^{[m_{i-1}+1,m_i]}$ as
$$V^{[m_{i-1}+1,m_i]}:=V^{\um_{i-1}+1}\oplus V^{\um_{i-1}+2}
\oplus\cdots\oplus V^{\um_i}.$$
Let $M(m_1,\ldots,m_r)$ denote the stabilizer of the decomposition
$$V=\bigoplus_{i=1}^{r}V^{[m_{i-1}+1,m_i]}.$$

\begin{lem} \label{intersection}
We have
\[M(m_1,\ldots,m_r)\cap U(\oA^{[m,1]})=\prod_{i=1}^rU(\oA^{[m_i,m_{i-1}]}).\]
\end{lem}
\begin{proof}
We set $e:=e(E|F)$. Let $l\in\{1,2,\ldots,m\}$ and
let $j\in\{0,1,\ldots,e-1\}$. Since the $\integers_E$-lattice chain
$\cL^l$ has period $1$, the equations~(\ref{di}) and (\ref{sum_di}) give
\[\dim_{k_E}L_j^l/L_{j+1}^l=\dim_{k_E}L_0^l/L_{1}^l
=\frac{N_l}{[E:F]}.\]
It follows that
\[d_j^l=\dim_{k_F}L_j^l/L_{j+1}^l=[k_E:k_F]\,\dim_{k_E}L_j^l/L_{j+1}^l=
[k_E:k_F]\,\frac{N_l}{[E:F]}=\frac{N_l}{e}.\]
Since $d_j^l=d_0^l=N_l/e$, we may and do fix an $\integers_F$-basis
$\cB^l:=(v^l_{0,1},\ldots,v^l_{0,N_l})$ of $\cL^l$, chosen to span $L^l_0$
over $\integers_F$. We put
\[v^l_{j,h}:=\begin{cases}
v^l_{0,h}&\text{ if $1\le h\le  \delta^l_j$,}\cr
\varpi_F\,v^l_{0,h}&\text{ if $\delta^l_j+1\le h\le N_l$,}
\end{cases}\]
where $$\delta_j^l\,:=\,\dim_{k_F}L^l_j/L^l_e\,=\,(e-j)d_j^l
\,=\,\frac{e-j}{e}\,N_l.$$
The $\integers_F$-lattice $L^l_j$ is then the $\integers_F$-linear span of
the set $\{v^l_{j,1},\ldots,v^l_{j,N_l}\}$, the cosets
$v_{j,h}^l+L_{j+1}^l$ ($1\le h\le N_l/e$) form a basis of
the $k_F$-space $L_j^l/L_{j+1}^l$, and
\[(v_{e-1,1}^l,\ldots,v_{e-1,N_l/e}^l, \ldots,
v_{1,1}^l,\ldots,v_{1,N_l/e}^l,v_{0,1}^l,\ldots,v_{0,N_l/e}^l)=\cB^l.\]
It follows that, for each $i\in\{1,2,\ldots,r\}$,
\[\cB^{[m_{i-1}+1,m_i]}:=
(\cB^{\um_{i-1}+1},\cB^{\um_{i-1}+2},\ldots,\cB^{\um_i})\]
is an $F$-basis of the vector space $V^{[m_{i-1}+1,m_i]}$ such that the cosets
$$v_{j,h}^{k+1}+L_{mj+k+1}^{[m,1]}, \quad\text{for $1\le h\le
N_{k+1}/e$},$$ form a
basis of the $k_F$-space
$$L_{mj+k}^{[m,1]}/L_{mj+k+1}^{[m,1]}\cong L_j^{k+1}/L_{j+1}^{k+1},$$ by
(\ref{goodquotient}).

Let $\cB$ denote the $F$-basis of $V$ defined as
\[\cB:=(\cB^{[1,m_1]},\cB^{[m_1+1,m_1+m_2]},\ldots,\cB^{[m_{r-1}+1,m_r]}).\]
We observe that we have by construction
\begin{equation} \label{samebasis}
\cB=\cB^{[1,m]},
\end{equation}
where $\cB^{[1,m]}$ is the $F$-basis corresponding to the partition $m$.

We will now use the basis $\cB$ to identify $A=\End_F(V)$ with $\Mat_N(F)$
and use the basis $\cB^{[\um_{i-1}+1,\um_i]}$ to identify
$\End_F(V^{[m_{i-1}+1,m_i]})$ with $\Mat_{N(i)}(F)$, where
\[N(i):=N_{\um_{i-1}+1}+N_{\um_{i-1}+2}+\cdots +N_{\um_i}.\]
Then $\oA^{[m,1]}$ becomes identified with the
matrices of the following form: the $(h,h')$-block has dimension
$$d_h(\cL^{[m,1]})\times
d_{h'}(\cL^{[m,1]}), \quad\text{if $0\le h,h'\le me-1$,}$$
and its entries lie in $\integers_F$ if $i\le i'$,
in $\fp_F$ otherwise.

Now the product $\prod_{i=1}^r\oA^{[\um_i,\um_{i-1}+1]}$ is viewed as diagonally
embedded in $\Mat_N(F)$, and $\oA^{[\um_i,\um_{i-1}+1]}$ becomes then
identified with the matrices of the following form: the
$(\um_{i-1}e+j,\um_{i-1}e+j')$-block has
dimension $$d_j(\cL^{[\um_i,\um_{i-1}+1]})\times
d_{j'}(\cL^{[\um_i,\um_{i-1}+1]}), \quad\text{if $0\le j,j'\le m_ie-1$,}$$
and its entries lie in $\integers_F$ if $j\le j'$, in $\fp_F$ otherwise.
Then the result follows from~(\ref{dLm1}).
\end{proof}

\section{Semisimple types} \label{section: semisimple}

Let $G=\GL(N,F)=\GL(V)$ and let $\fs=[M,\sigma]_G$ be a point in the Bernstein
spectrum $\mathfrak{B}(G)$. The Levi subgroup $M$ is the stabilizer of a
decomposition
$V=\bigoplus_{l=1}^{m}V^l$ of $V$ as a direct sum of nonzero
subspaces $V^l$. We set $N_l:=\dim_F V^l$, and
$A_l:=\End_F(V^l)\cong\Mat_{N_l}(F)$.
Then $N_1+\cdots+N_m=N$, and $M$ is isomorphic to
$\GL(N_1,F)\times\cdots\times\GL(N_m,F)$, and the
supercuspidal representation $\sigma$ of $M$ is of the form
$\sigma=\pi_1\otimes\cdots\otimes\pi_m$, where
$\pi_l$ is an irreducible supercuspidal representation of the group $\GL(N_l,F)$,
for $l=1,\ldots,m$.
We set $\ft:=[M,\sigma]_M$.

By \cite[Theorem~(8.4.1)]{BK3}, for each $l$, there is a maximal simple type
$(J^l,\lambda^l)$ occuring in $\pi_l$. The pair
$(J_M,\tau_M):=(J^1\times \cdots \times J^m,\lambda^1\otimes\cdots
\otimes\lambda^m)$ is then an $\ft$-type in $M$.

By definition (see \cite[(5.5.10)]{BK3}), for each $l$, there
exists an element $\beta_l\in A^l$ for which the algebra $E_l:=F[\beta_l]$ is
a field and a principal $\integers_F$-order $\oA^l$ in $A^l$, of period
$e(E_l|F)$, with Jacobson radical $\oP_l$, such that
\[J^l=\begin{cases}
J(\beta_l,\oA^l)&\text{(as defined in \cite[(3.1.14)]{BK3}) if
$\beta_l\notin F$,}\cr
U(\oA^l) &\text{if $\beta_l\in F$.}
\end{cases}\]
For each $x\in A^l$, we will write
\begin{equation} \label{valuation}
\nu_{\oA^l}(x):=\max\left\{n\in\Zset\,:\,x\in \oP_l\right\}.
\end{equation}
Let $\cL^l$ denote the $\integers_E$-lattice chain defining the maximal
$\integers_E$-order $\oB^l:=\oA^l\cap\End_E(V^l)$.
We have
\begin{equation} \label{Jquotient}
J(\beta,\oA^l)/J^1(\beta,\oA^l)=U(\oB^l)/U^1(\oB^l)\cong\GL(f_l,k_E).
\end{equation}

\subsection{Simple types} \label{subsection: simple}

We assume in this subsection that the $N_l$ are all equal to $N/m$ and that
$\pi_l\cong\pi_j\chi_j$, with
$\chi_j$ an unramified character of $\GL(N/m,F)$, for each
$l,j\in\{1,\ldots,m\}$.
In particular, $M$ is then isomorphic to $\GL(N/m,F)^m$, and
by \cite[Theorem~(8.4.2)]{BK3}, we can assume that all the $\beta_l$,
all the $\oA^l$, all the $\cL^l$, all the $J^l$ and all the $\lambda^l$
are equal. We will denote by $E$ (resp. $\beta$) the common value of the $E_l$
(resp. $\beta_l$).

Using the second addition procedure~\ref{second}, we define:
the $\integers_E$-lattice chain
\begin{equation} \label{la-chaine}
\cL:=\cL^m+\cL^{m-1}+\cdots+\cL^1,\;\;\text{ and }\;\;
\oA:=\End_{\integers_F}^0(\cL).\end{equation}
If $\beta$ belongs to $F$, we set $J:=U(\oA)$. Otherwise, let
$n:=-\nu_{\oA^1}(\beta)$,
then $[\oA,mn, 0,\beta]$ is a simple stratum in the sense of
\cite[Definition~(1.5.5)]{BK3}, let $(J,\lambda):=(J(\beta,\oA),\lambda)$
be the corresponding simple type in $G$.

Let $\oB$ denote the principal $\integers_E$-order in
$B:=\Mat_{N/[E:F]}(E)$ defined by $\oB:=B\cap\oA$.
We have $m=e(\oB)=e(\oB|\integers_E)$.
In the case when $\beta\in F$, we have $m=e(\oA)$.

\begin{defn} \label{Jssimple} We set
\[\oA^\fs:=\oA(N/m,\ldots,N/m)\;\;\text{ and }\;\;
J^\fs:=U(\oA^\fs),\]
where $\oA(N/m,\ldots,N/m)$ is defined by~Definition~\ref{Ordre}.
\end{defn}

\begin{lem} \label{ordre_inclusI}
The $\integers_F$-order $\oA$ is contained in the $\integers_F$-order
$\oA^\fs$.
\end{lem}
\begin{proof}
We have $\oA=\oA(N/e(\oA),\ldots,N/e(\oA))$.
In the case when $J=U(\oA)$, we have $\oA^\fs=\oA$. Otherwise, the
statement follows immediately from the above descriptions of the
orders $\oA$, $\oA^\fs$, and from the fact (see \cite[Proposition~(1.2.4)]{BK3})
that
\[e(\oA)=m \cdot e(E|F).\]
Indeed, from the above descriptions of the orders $\oA$, $\oA^\fs$, we have
\begin{equation} \label{unipotents}
\oA^\fs\cap U=\oA\cap U,\;\;\oA^\fs\cap \overline U=\oA\cap \overline
U,
\end{equation}
\begin{equation} \label{Levi}
M\cap\oA^\fs\,\cong\,(\GL(N/m,\integers_F))^m,
\end{equation}
while $M\cap\oA$ is isomorphic to the product of $m$ copies of the order of
$e(E|F)\times e(E|F)$ blocks matrices of the following form: the $(j,l)$-block
has dimension $N/e(\oA)\times N/e(\oA)=(N/e(E|F)m\times N/e(E|F)m)$,
$0\le j,l\le e(E|F)-1$, and its entries lie in
$\integers_F$ if $j\le l$, in $\varpi_F\integers_F$ otherwise, so that
$M\cap\oA\subset M\cap\oA^\fs$.
\end{proof}

We set
\begin{equation}
f=\frac{N}{[E:F]\cdot m}.
\end{equation}
Let $K/E$ be an unramified field extension of degree $f$ with
\[K^\times\,\subset\,\fK(\toB),\]
where $\fK(\toB)$ is defined by~(\ref{normalizer}),
and let $C=\End_K(V)\cong\Mat_{m}(K)$.
We view $\varpi_E$ as a prime element of $K$.
For $i=1,\ldots, m-1$, let $s_{i,C}$ denote the matrix in $C$ of the
transposition $i\leftrightarrow i+1$, that is,
\[s_{i,C}=
\left(\begin{matrix}
\Id_{i-1}&&&\cr
&0&1&\cr
&1&0&\cr
&&&\Id_{m-i-1}
\end{matrix}\right),\]
and let $s_{0,C}=\Pi_{m,C}s_{1,C}\Pi_{m,C}^{-1}$, with
\[\Pi_{m,C}=\left(\begin{matrix}
0&\Id_{m-1}\cr
\varpi_E&0\end{matrix}\right).\]
We fix the embedding
\[\bigotimes\Id_{N/m}\colon C\hookrightarrow \Mat_N(K)\quad
c=(c_{ij})\mapsto c\otimes \Id_{N/m}=(c_{ij}\Id_{N/m}),\]
$c\otimes \Id_{N/m}$ being a block matrix with scalar blocks.

Let $\widetilde W_C$ be the group generated by
\[S=\{s_{0,C}\otimes \Id_{N/m},s_{1,C}\otimes \Id_{N/m},\ldots,
s_{m-1,C}\otimes \Id_{N/m}\}.
\]
Then $(\widetilde W_C,S)$ is a Coxeter group of type $\tilde A_{m-1}$.
\begin{thm} \label{simplet}
The representation $\alpha=\Ind_J^{J^\fs}(\lambda)$ is irreducible.
Hence the pair $(J^\fs,\alpha)$ is an $\fs$-type.
\end{thm}
\begin{proof}
In the case when $J=U(\oA)$, we have $J^\fs=J$, so the result follows
trivially in this case. We will assume from now on that $J=J(\beta,\oA)$.
For any $i\in\{1,\ldots,m-1\}$,
\[s_{i,C}\otimes \Id_{N/m}=\left(\begin{matrix}
\Id_{(i-1)N/m}&&&\cr
&0&\Id_{N/m}&\cr
&\Id_{N/m}&0&\cr
&&&\Id_{(m-i-1)N/m}
\end{matrix}\right)\,\notin\,J^\fs,\]
and
\[\Pi_{m,C}\otimes \Id_{N/m}=\left(\begin{matrix}
0&\Id_{(m-1)N/m}\cr
\varpi_E\Id_{N/m}&0\end{matrix}\right)\,\notin\,J^\fs).\]
Hence ${\widetilde W}_C\cap J^\fs=\{1\}$, which gives
\begin{equation} \label{entrelacementWB}
J^\fs\,\cap\,(J\cdot {\widetilde W}_C\cdot J)=J.
\end{equation}
Then the result follows from the fact (see
\cite[Propositions~(5.5.11) and (5.5.14)~(iii)]{BK3}) that
\[I_G(\lambda)\,\subset\,J\cdot {\widetilde W}_C\cdot J.\]
\end{proof}

\subsection{In the Levi subgroup $\tM$} \label{subsection: tM}

We will now consider the case of an arbitrary point $\fs=[M,\sigma]_G$ in
$\fB(G)$, with $G=\GL(N,F)$.
Let $\tM$ denote the unique Levi subgroup of $G$ which contains
$\Nor_\ft$ (see~\ref{Nor_t}) and is minimal for this property.

We write $\sigma=\pi_1\otimes\pi_2\otimes\cdots\otimes\pi_m$ as
$$\sigma=(\sigma_1, \ldots, \sigma_1, \sigma_2,\ldots,\sigma_2,
\ldots, \sigma_t, \ldots,\sigma_t),$$ where $\sigma_j$, a
supercuspidal representation of $\GL(N'_j,F)$, is repeated $\e_j$
times, $1 \leq j \leq t$, and $\sigma_1$, $\ldots$, $\sigma_t$ are
pairwise distinct (after unramified twist). The integers $\e_1$,
$\e_2$, $\ldots$, $\e_t$ are called the \emph{exponents} of
$\sigma$. Then we have
\[M\,\cong\,\GL(N_1',F)^{\e_1}\times\GL(N'_2,F)^{\e_2}\times\cdots\times
\GL(N_t',F)^{\e_t},\]
and
\[\tM\,\cong\,
\GL(\e_1N'_1,F)\times\GL(\e_2N_2',F)\times\cdots\times\GL(\e_tN_t',F).\]

For every $j\in\{1,\ldots,t\}$, we set
\[\fs_j=[\GL(N_j',F)^{\e_j},\sigma_j^{\otimes \e_j}]_{\GL(\e_jN'_j,F)}.\]
Then let $(K^j,\tau^j)$ be
the $\fs_j$-type in $\GL(\e_jN'_j,F)$ (a simple type) defined as
in the previous section, and
let $(\tK^j,\ttau^j)$ be the ``modified simple type'' attached to
$(K^j,\tau^j)$ as in \cite[proof of Prop.~1.4]{BK2}.

\begin{lem} \label{modifiedst}
We have $\tK^j\subset J^{\fs_j}$ and $\alpha_i=
\Ind_{\tK^j}^{J^{\fs_j}}(\ttau^j)$ is irreducible.
\end{lem}
\begin{proof}
There is an isomorphism of Hecke algebras
\[\cH(\GL(\e_jN'_j),\ttau^j)\,\cong
\,\cH(\GL(\e_jN_j),\tau^j)\]
such that, if $\tilde f\in \cH(\GL(\e_jN'_j),\ttau^j)$ has support
$\tK^j g \tK^j$, for some element $g\in\GL(\e_jN_j',F)$, then its image
$f$ in $\cH(\GL(\e_jN'_j),\tau^j)$ has support
$K^jgK^j$ (see \cite[(7.2.19)]{BK3}). Then the result follows
from Theorem~\ref{simplet}.
\end{proof}

We set
\[\fs_{\tM}=[M,\sigma]_{\tM},\;\;\;\;
\widetilde\oA^{\fs}=\oA^{\fs_1}\times\cdots\times\oA^{\fs_t},\;\;\;\;
\tJ^{\fs}=U(\widetilde\oA^{\fs}),\]
\begin{equation} \label{tK}
\tK=\tK^1\times\cdots\times\tK^t\,\subset\,\tM,\;\;\;\;
\ttau=\ttau^1\otimes\cdots\otimes\ttau^t.\end{equation}
Note that
\begin{equation} \label{tMtJ}
\tJ^{\fs}=\tM\cap J^\fs.
\end{equation}

It immediately follows from Lemma~\ref{modifiedst} that:
\begin{lem} \label{modifiedstM}
We have $\tK\subset \tJ^{\fs}$ and
$\talpha=\Ind_{\tK}^{\tJ^{\fs}}(\ttau)$ is irreducible.
\end{lem}


\subsection{Review of endo-classes}
\label{subsection: endoclass}

We recall that a \emph{simple pair $(k,\beta)$ over $F$} consists of an
integer $k$ and a nonzero element $\beta$ generating a field extension $E$
of $F$ such that
\[-k>\max\left\{k_0(\beta,\oA(E)),\nu_E(\beta)\right\},\]
where $\nu_E$ is the standard additive valuation on $E$ and
$k_0(\beta,\oA(E))$ is defined by \cite[(1.4.5)]{BK3}, with
$\oA(E)$ denoting the unique hereditary $\integers_F$-order in $\End_F(E)$
such that $\fK(\oA(E))\supset E^\times$.

Let $(k,\beta)$ be a given simple pair in which $k\ge 0$.
A \emph{ps-character} (attached to the simple pair $(k,\beta)$) is then a
triple $(\Theta,k,\beta)$, where $\Theta$ is a simple-character-valued
function, such that to each triple $(V,\oB,m)$, where $V$ is a
finite-dimensional $E$-vector space, $\oB$ is a hereditary
$\integers_E$-order in $\End_E(V)$, and $m$ is an integer such that
$[m/e(\oB|\integers_E)]=k$, the function $\Theta$ attaches a simple
character $\Theta(\oA)\in\cC(\oA,m,\beta)$, called the \emph{realization
of $\Theta$ on $\oA$ of order $m$}. (If we put $n:=-\nu_E(\beta)\,e(\oB)$,
the stratum $[\oA,n,m,\beta]$ is simple and the simple character set
$\cC(\oA,m,\beta)$ of \cite[(3.2)]{BK3} is defined.)

These realizations are subject to the
following coherence condition: if we have two realizations $\Theta(\oA_1)$
and $\Theta(\oA_2)$ of on orders $\oA_1$, $\oA_2$, they are related by
$\Theta(\oA_2)=\btau_{\oA_1,\oA_2,\beta}(\Theta(\oA_1))$, where
\[\btau_{\oA_1,\oA_2,\beta}\colon\cC(\oA_1,m,\beta)\to\cC(\oA_2,m,\beta)\]
is the canonical bijection of \cite[(3.6.14)]{BK3}.

Following \cite[\S 4.3]{BK2}, we will say that two ps-characters
$(\Theta_1,k_1,\beta_1)$ and $(\Theta_2,k_2,\beta_2)$ are
\emph{endo-equivalent} if there exists an $F$-vector space $V$, hereditary
$\integers_F$-orders $\oA_1$, $\oA_2$ in $\End_F(V)$, and realizations
$\Theta_i(\oA_i)$ of the $\Theta_i$ of same level, such that $\oA_1\cong
\oA_2$ as $\integers_F$-orders, and such that the simple characters
$\Theta_i(\oA_i)$ intertwine in $\Aut_F(V)$. Endo-equivalence in
equivalence relation on the set of ps-characters over $F$. One refers to
the equivalence classes as \emph{endo-classes} of simple characters.

If the supercuspidal representation $\pi_l$ of $\GL(N_l,F)$ contains the
trivial character of $U^1(\oA^l)=1+\oP_l$,
then $\pi_l$ is said to be of level-zero. Otherwise,
there exists a simple stratum
$[\oA^l,n_l,0,\beta_l]$ in $A_l$ and a
simple character $\theta_l\in\cC(\oA^l,0,\beta_l)$ such that the
restriction of $\lambda^l$ to $H^1(\beta_l,\oA^l)$ is a multiple of
$\theta_l$. (Here $H^1(\beta_l,\oA^l)$ is defined as in
\cite[(3.1.14)]{BK3}.) Since $[\oA^l,n_l,0,\beta_l]$ is simple, we have
$n_l=-\nu_{\oA^l}(\beta_l)$. Then
each representation $\lambda^l$ is given as follows.
There is a unique irreducible representation $\eta_l$ of $J^1(\beta,\oA^l)$
whose restriction to $H^1(\beta,\oA^l)$ is a multiple of $\theta_l$. The
representation $\eta_l$ extends to a representation $\kappa_l$ which is a
$\beta$-extension of $\eta_l$, and we have
$\lambda^l=\kappa_l\otimes\rho_l$, where $\rho_l$ is the inflation of an
irreducible representation of $\GL(f_l,k_E)$, with $f_l$ defined
by~(\ref{Jquotient}).

If the representation $\pi_l$ is of level zero, we set
$\Theta_{\pi_l}=\{\Theta^0\}$, where
$\Theta^0$ is the trivial ps-character (that is, if $\oA$ is a
hereditary $\integers_F$-order in some $\End_F(V)$, the realization of
$\Theta^0$ on $\oA$ is the trivial character of $U^1(\oA)$).
Otherwise, the simple character $\theta_i$ determines a ps-character
$(\Theta_l,0,\beta)$ and hence an endo-class $\Theta_{\pi_l}$.

We will denote by $\Theta(1)$, $\Theta(2)$, $\ldots$, $\Theta(q)$ the distinct
endo-classes arising in the set $\{\Theta_{\pi_1},\ldots,\Theta_{\pi_m}\}$.

\subsection{The homogeneous case} \label{homogeneous_case}

In this subsection, we assume that all the representations $\pi_1$,
$\pi_2$, $\ldots$, $\pi_m$ admit the same endo-class.
It follows that all the elements $\beta_1$, $\ldots$, $\beta_m$ may be assumed
to be equal.
We will denote by $E$ (resp. $\beta$) the common value of the $E_l$
(resp. $\beta_l$).

Let $l\in\{1,\ldots,m\}$, and let $(v_1^l,v_2^l,\ldots,v_{N_l}^l)$ be an
$F$-basis of $V^l$, with respect to which $\oA^l=\oA(\cL^l)$ is identified
with $\oA(N_l/e(E|F),\ldots,N_l/e(E|F))$.
We have $L^l_0=\integers_Fv_1^l\oplus\cdots\oplus\integers_Fv_{N_l}^l$.
We set \[L_{i,\max}^l:=\fp^i L_0^l,\;\;\text{for any $i\in\Zset$.}\]
Then
\begin{equation} \label{Lmax1}
\cL^l_{\max}:=\left\{ L_{i,\max}^l\,:\,i\in\Zset\right\}
\end{equation}
is an $\integers_F$-lattice chain in $V^l$ of period $1$, and we have
\[\oA(\cL^l_{\max}):=\End_{\integers_F}^0(\cL^l_{\max})=\oA(N_l)=
\Mat_{N_l}(\integers_F)\,\supset\,\oA^l.\]

Following the second addition procedure defined in the
subsection~\ref{addition_procedure}, we assemble the $\integers_E$-lattices
chains $\cL^1$, $\ldots$, $\cL^m$ into the $\integers_E$-lattice chain
\begin{equation} \label{Lchaine}
\bar\cL:=\cL^m+\cL^{m-1}+\cdots+\cL^1
\end{equation} in $V$, of period $m$, and
we assemble the $\integers_F$-lattices chains $\cL^1_{\max}$, $\ldots$,
$\cL^m_{\max}$ into the $\integers_F$-lattice chain
\begin{equation} \label{Lmax}
\bar\cL_{\max}:=\cL^m_{\max}+\cL^{m-1}_{\max}+\cdots+\cL^1_{\max}
\,=\,\left\{
\bar L_{\max,i}\,:\,i\in\Zset\right\}
\end{equation} in $V$, of period $m$.
Let $j\in\Zset$ and $k\in\{0,1,\ldots,m-1\}$.
From~(\ref{goodquotient}), we have
\[\bar L_{\max,mj+k}/\bar L_{\max,mj+k+1}\,
\cong \,L_{\max,j}^{k+1}/L_{\max,j+1}^{k+1}.
\]
Hence:
\begin{equation} \label{dimension_Lmax}
d_{mj+k}(\bar\cL_{\max})=N_{k+1}.\end{equation}
It follows that
\[\oA(\bar\cL_{\max})=\oA(N_1,N_2,\ldots,N_m).\]
We put
\begin{equation} \label{BB}
B:=\End_E(V)\;\;\text{ and }\;\;\oB:=\End_{\integers_E}^0(\bar\cL).
\end{equation}
Considering $\bar\cL$ as an $\integers_F$-lattice chain, we put
\begin{equation} \label{AA}
\oA:=\End_{\integers_F}^0(\bar\cL).\end{equation}
We have $\oB=\oA\cap B$.

\smallskip

The following definition, lemma and theorem generalize
Definition~\ref{Jssimple}, Lemma~\ref{ordre_inclusI}, and
Theorem~\ref{simplet}, respectively.

\begin{defn} \label{Jshomogene}
We set
\[\oA^{\fs}:=\oA(N_1,N_2,\ldots,N_m)\;\;\text{ and }\;\;
J^\fs:=U(\oA^\fs).\]
\end{defn}

\begin{lem} \label{ordre_inclusII}
The $\integers_F$-order $\oA$ is contained in the $\integers_F$-order $\oA^{\fs}$.
\end{lem}
\begin{proof}
We have
\[\oA^\fs\cap U=\oA\cap U,\;\;\oA^\fs\cap \overline U=\oA\cap \overline U,
\]
\[
M\cap\oA^\fs\,\cong\,\prod_{l=1}^m\GL(N_l,\integers_F).
\]
In the notation of subsection~\ref{second}, setting
$(m_1,\ldots,m_r)=(1,\ldots,1)$, we have $r=m$, $M=M(m_1,\cdots,m_r)$.
Then $\um_{l-1}=l-1=\um_l-1$,
hence $\bar\cL^{[\um_l,\um_{l-1}+1]}=\cL^l$,
\[\oA^l=\oA(\cL^l)\cong\oA(N_l/e(E|F),\ldots,N_l/e(E|F)),\]
and Lemma~\ref{intersection} gives
\[M\cap\oA\cong
\prod_{l=1}^m U(\oA^l).\]
Since $U(\oA^l)\subset \GL(N_l,\integers_F)$, the result follows.
\end{proof}

\smallskip

We set
\begin{equation} \label{defn_n}
n:=\max(n_1,\ldots,n_m).
\end{equation}

\begin{lem} \label{strate_simple}
With notation as above, $[\oA,nm,0,\beta]$ is a simple stratum.
\end{lem}
\begin{proof}
We have to check that the four conditions occuring in
\cite[Definition~(1.5.5)]{BK3} are satisfied.
\begin{itemize}
\item[(i)]
We know that the algebra $E=F[\beta]$ is a field, since the strata
$[\oA^l,n_l,0,\beta]$ are simple.
\item[(ii)]
We defined $\bar\cL=\{\bar L_i\,:\,i\in\Zset\}$ to be an $\integers_E$-lattice
chain in the $E$-vector space $V$.
Hence, by \cite[Proposition~(1.2.1)]{BK3}, we have $E^\times\subset\fK(\oA)$.
\item[(iii)]
Let $l\in\{1,\ldots,m\}$.
We set $\oQ_l:=\oB_l\cap\oP_l$.
Since $\nu_{\oA^l}(\beta)=-n_l$, the definition~(\ref{valuation}) for
$\nu_{\oA^l}$ shows that
\[\beta\in\End_E(V^l)\cap\oP_l^{-n_l}\;\;\text{ and }\;\;
\beta\notin\oP_l^{-n_l+1},\]
that is,
\[\beta\in\oQ_l^{-n_l}\;\;\text{ and }\;\;
\beta\notin\oQ_l^{-n_l+1}.\]
By \cite[Proposition~(1.2.4)]{BK3}, we know that $\oQ_l$ is the Jacobson radical
of the $\integers_E$-order $\oB_l$.
Hence $\oQ_l^i=\End_{\integers_E}^i(\cL^l)$ for each $i\in\Zset$, and
$\beta(L_j^l)$ is contained in $L^l_{j-n_l}$ and not in
$L^l_{j-n_l+1}$.
Now, it follows from~(\ref{goodchain}) that
\[\beta(L_{mj+k})=
\beta(L_{j+1}^1)\oplus\cdots\oplus \beta(L_{j+1}^{k})\oplus
\beta(L_{j}^{k+1})\oplus\cdots\oplus \beta(L_{j}^m),\]
for each $j\in\Zset$ and each $k\in\{0,1,\ldots,m-1\}$.
Since $n=\max(n_1,\ldots,n_m)$, we have $L_{j-n_l}^l\subset L^l_{j-n}$,
for each $l$. It gives $\beta(L_{mj+k})\subset L_{m(j-n)+k}$.
On the other side there exists $l_0\in\{1,\ldots,m\}$ such that $n=n_{l_0}$, and
hence $\beta(L_j^{l_0})$ is not contained in $L^l_{j-n+1}$.
It follows that $\beta(L_{mj+k})$ is not contained in $L_{m(j-n)+k+1}$, that is,
\[\beta\in\oQ^{-n}\;\;\text{ and }\;\;
\beta\notin\oQ^{-n+1},\]
where $\oQ$ denotes the Jacobson radical of $\oB$. Since
$\oQ^i=B\cap\oP^i$ for each $i\in\Zset$ (by
\cite[Proposition~(1.2.4)]{BK3}), we get
$\nu_\oA(\beta)=-nm$.
\item[(iv)]
Let $A(E):=\End_F(E)$. The algebra $A(E)$ contains the principal
$\integers_F$-order
\[\oA(E):=\End_{\integers_F}^0
\left(\left\{\fp_E^i\,:\,i\in\Zset\right\}\right).\]
We have $E^\times\subset\fK(\oA(E))$ and
\cite[Proposition~(1.4.13)~(ii)]{BK3} gives
\[k_0(\beta,\oA)=mk_0(\beta,\oA(E))=k_0(\beta,\oA^l),\;\;\text{ for each
$l\in\{1,\ldots,m\}$.}\]
Since $[\oA^l,n_l,0,\beta]$ is a simple stratum, we have
\[0<-k_0(\beta,\oA^l).\] Hence $0<-k_0(\beta,\oA)$ and
$[\oA,mn,0,\beta]$ is simple.
\end{itemize}
\end{proof}

Since $[\oA,nm,0,\beta]$ is a simple stratum, we can associate to it
the compact open subgroups $J(\beta,\oA)$ and $H^1(\beta,\oA)$ of
$U(\oA)$, defined following \cite[(3.1.14)]{BK3}.

As in \cite[\S 7.1, 7.2]{BK3}, the set
\begin{equation} \label{K}
K:=H^1(\beta,\oA)\cap \overline U\,\cdot\, J(\beta,\oA)\cap P
\end{equation}
is then a subgroup of $U(\oA)$ containing $H^1(\beta,\oA)$.

Definition~\ref{Jshomogene} and Lemma~\ref{ordre_inclusII} imply
\begin{equation} \label{JJ}
K\subset J^{\fs}.
\end{equation}
As in \cite[\S 7.2.1]{BK2}, it admits an irreducible representation $\kappa$,
trivial on $K\cap\overline U$, $K\cap U$,
whose restriction to $H^1(\beta,\oA)$ is a multiple of
$\theta=\Theta(\oA)$, and such that $\kappa_{|K\cap M}$ is of the form
$\kappa_1'\otimes\cdots\otimes\kappa_m'$ for some $\beta$-extension
$\kappa'_l$ of $\eta_l$.
As in \cite[\S 7.2]{BK2}, we can choose the decomposition
$\lambda^l=\kappa_l\otimes\rho_l$ above so that $\kappa_l=\kappa_l'$; we
assume this has been done. We have canonically
\[K/K\cap J^1(\beta,\oA)\,\cong\,\prod_{l=1}^m
J(\beta,\oA^l)/J^1(\beta,\oA^l)\,\cong\,\prod_{l=1}^m\GL(f_l,k_E),\]
and we can inflate the cuspidal representation
$\rho_1\otimes\cdots\otimes\rho_m$ of $\prod_{l=1}^m\GL(f_l,k_E)$ to a
representation $\rho$ of $K$ and form
\begin{equation} \label{tauGL}
\tau=\kappa\otimes\rho.
\end{equation}

Moreover similar proofs of those of \cite[Theorem~7.2.1, Main Theorem~8.2]{BK2}
show that $(K,\tau)$ is a $G$-cover of the pair $(\tK,\ttau)$ defined
in (\ref{tK}) and give the following formula for the intertwining:
\begin{equation} \label{GI}
\cI_G(\tau)=K\cdot\cI_{\tM}(\ttau)\cdot K.
\end{equation}

\begin{thm} \label{homogenet}
Let $J^\fs$ be as in Definition~\ref{Jshomogene}.
Then the representation $\alpha:=\Ind_K^{J^\fs}(\tau)$ is irreducible.
Hence the pair $(J^\fs,\alpha)$ is an $\fs$-type.
\end{thm}
\begin{proof}
Using equations~(\ref{JJ}) and~(\ref{GI}), we obtain
\[\cI_{J^\fs}(\tau)=K\cdot\cI_{\tM\cap J^{\fs}}(\ttau)\cdot K.\]
On the other side, equation~(\ref{tMtJ}) and Lemma~\ref{modifiedstM} imply that
\[\cI_{\tM\cap J^{\fs}}(\ttau)=\tK\,\subset K.\]
Hence $\cI_{J^\fs}(\tau)=K$, and the result follows from
Proposition~\ref{IK}.
\end{proof}

\subsection{The general case} \label{general_case}

The Levi subgroup $\tM$ defined in the beginning of the
subsection~\ref{subsection: tM} is the $G$-stabilizer of a decomposition
\[V=\tV^1\oplus\tV^2\oplus+\cdots\oplus\tV^t,\]
of $V$ as a direct sum of nonzero subspaces $\tV^j$.

Since the endo-class of a supercuspidal representation only depends on the
corresponding point in the Bernstein spectrum (see \cite[Proposition~4.5]{BK2}),
we can associate to each $\tV^j$ an endo-class of simple characters,
namely $\Theta_{\pi_l}$ for any $l$ such that $V^l\subset\tV^j$.

Now let $\barM\supset M$ be the Levi subgroup in $G$ defined as in
\cite[\S 8.1]{BK2}, that is, for each $i$, let $\bar V^i$ be the sum of
those $\tV^j$ whose associate endo-class $\Theta_{\pi_j}$ is $\Theta(i)$,
and write $\barM$ for the $G$-stabilizer of a decomposition
\[V=\bar V^1\oplus\bar V^2\oplus\cdots\oplus\bar V^q.\]
Setting $\barN_i:=\dim_F\bar V^i$, we get
$$\barM\cong\GL(\barN_1,F)\times\cdots\times\GL(\barN_q,F).$$
We put
\begin{equation} \label{bK}
\barK:=K_{1}\times K_{2}\times\cdots\times K_{q}\;\;\text{ and }\;\;
\bartau:=\tau_{1}\times \tau_{2}\times\cdots\times \tau_{q},
\end{equation}
where the pairs $(K_{i},\tau_{i})$ are defined as in~(\ref{K}), (\ref{tauGL}).
Then a similar proof as those of \cite[\S 7.2]{BK2} shows that the pair
$(\barK,\bartau)$ is a $\barM$-cover of $(J_M,\tau_M)$.

\smallskip

For each $i\in\{1,\ldots,q\}$, let $\bar\cL^i$, $\bar\cL_{\max}^i$
respectively denote the $\integers_E$-lattice chain in the $E$-vector
space $\bar V^i$ defined
by~(\ref{Lchaine}), and the $\integers_F$-lattice chain in the $F$-vector
space $\bar V^i$
defined by~(\ref{Lmax}). Let $m_i$ denote the number of representations
$\pi_l$ ($1\le l\le m$) with endo-class $\theta_i$. Then $\bar\cL^i$,
considered as an $\integers_F$-lattice chain, has period
$e_i:=e(\bar\cL^i)=e(E_i|F)\,m_i$, and $\bar\cL_{\max}^i$
has period $e(\bar\cL^i_{\max})=m_i$.

Then let $\Lambda^i$ (resp. $\Lambda^i_{\max}$) denote the (strict) \emph{lattice
sequence} defined by the lattice chain $\bar\cL^i$ (resp. $\bar\cL_{\max}^i$),
considered as $\integers_F$-lattice chains.
Then, using the addition of lattice sequences recalled in~(\ref{seq_add}), we
define
\begin{equation} \label{Lambda}
\Lambda:=\Lambda^1\,\oplus\,\Lambda^2\,\oplus\,\cdots\,
\oplus\,\Lambda^q,\end{equation}
and
\begin{equation} \label{Lambda_max}
\Lambda_{\max}:=e(E_1|F)\Lambda_{\max}^1\,\oplus\,
e(E_2|F)\Lambda_{\max}^2\,\oplus\,\cdots
\,\oplus\,e(E_q|F)\Lambda^q_{\max}.
\end{equation}

Let $\cL_{\Lambda}$,
$\cL_{\Lambda_{\max}}$ denote the $\integers_F$-lattice chains attached to
the lattice sequences $\Lambda$, $\Lambda_{\max}$, respectively, as
in~(\ref{LL}).
Let $\oA_\Lambda$, $\oA_{\Lambda_{\max}}$ denote the hereditary
$\integers_F$-orders in $A$ defined by the lattice chain $\cL_\Lambda$,
$\cL_{\Lambda_{\max}}$, respectively. We have (see
\cite[Proposition~2.3.~(i)]{BK2}):
\begin{equation} \label{orderLambda}
\oA_\Lambda=\oA(\cL_{\Lambda})=\fa_0(\Lambda)\;\;\text{ and }\;\;
\oA_{\Lambda_{\max}}=\oA(\cL_{\Lambda_{\max}})=\fa_0(\Lambda_{\max}),
\end{equation}
where $\fa_0(\Lambda)$, $\fa_0(\Lambda_{\max})$ are defined as
in~(\ref{aLambda}).

\begin{lem} \label{inclusionLL}
The $\integers_F$-order $\oA_\Lambda$ is contained in the
$\integers_F$-order $\oA_{\Lambda_{\max}}$.
\end{lem}
\begin{proof}
Let $e:=\lcm\{e_1,\ldots,e_q\}$. Both $\Lambda$ and $\Lambda_{\max}$ have
period $e$.
>From (\ref{seq_add}), we have
\[\Lambda(ex)=\Lambda^1(e_1x)\oplus\cdots\oplus\Lambda^q(e_qx),\;\;
\Lambda_{\max}(ex)=\Lambda^1_{\max}(e_1x)\oplus\cdots\oplus
\Lambda^q_{\max}(e_qx),\]
for each $x\in\Rset$. On the other side, (\ref{reel}) gives
\[\Lambda^i\left(\frac{e_i}{e}j\right)=\Lambda(l_i(j))\;\;\text{ and }\;\;
\Lambda^i_{\max}\left(\frac{e_i}{e}j\right)=\Lambda_{\max}(l_i(j)),\]
for each $j\in\Zset$,
where $l_i(j)$ is the integer defined by the relation
$$l_i(j)-1<\displaystyle\frac{e_i}{e}j\le l_i(j).$$
Hence
\[\Lambda(j)=\Lambda^1(l_1(j))\oplus\cdots\oplus\Lambda^q(l_q(j)),\;\;
\Lambda_{\max}(j)=\Lambda^1_{\max}(l_1(j))
\oplus\cdots\oplus\Lambda^q_{\max}(l_q(j)).\]
Then the result is consequence of Lemma~\ref{ordre_inclusII}.
\end{proof}

The following definition generalizes Definitions~\ref{Jssimple}
and~\ref{Jshomogene}.
\begin{defn} \label{Jsgeneral}
We set
\[\oA^\fs:=\oA_{\Lambda_{\max}},\;\;\text{ and }\;\;
J^{\fs}:=U(\oA^\fs).\]
\end{defn}

\begin{example} \label{Ex2}
We assume here that $q=m$, that is, the representations $\pi_l$ have all
distinct endo-classes. It implies that $\barM=\tM=M$. Then each lattice
sequence $\Lambda_{\max}^l$ has period $1$, and so $\Lambda_{\max}$ has
also period $1$. We get in this case $J^{\fs}=\GL(N,\integers_F)$.
\end{example}

\begin{thm} \label{ssinclusion}
There exists a $G$-cover $(J,\tau)$ of $(J_M,\lambda_M)$ such that
\begin{itemize}
\item $J\subset J^\fs$,
\item
$\alpha:=\Ind_J^{J^\fs}(\tau)$ is irreducible.
Hence $(J^\fs,\alpha)$ is an $\fs$-type.
\end{itemize}
\end{thm}
\begin{proof}
Let $(J,\tau)$ be the $G$-cover of $(\barK,\bartau)$ constructed in the
similar way as in \cite[\S 8]{BK2}, in particular, we have
\[J\subset U(\oA_\Lambda).\]
Then the first assertion follows from Lemma~\ref{inclusionLL}.

On the other side the same proof as those of \cite[\S 8.2, Main Theorem]{BK2}
gives the following formula for the intertwining:
\[\cI_G(\tau)=J\cdot \cI_{\tM}(\tau_{\tM})\cdot J.\]
Since $J\subset J^\fs$, it implies:
\[\cI_{J^\fs}(\tau)=J\cdot \cI_{\tM\cap J^\fs}(\tau_{\tM})\cdot J=
J\cdot \cI_{\tJ^\fs}(\ttau)\cdot J.\]
Now, by Lemma~\ref{modifiedstM}, we have
\[\cI_{\tJ^\fs}(\tau_{\tM})=\tK.\]
We get
\[\cI_{J^\fs}(\tau)=J,\]
and the result follows from Proposition~\ref{IK}.
\end{proof}

\section{Supercuspidal Bernstein components}

Let $\fs=[G,\pi]_G$, where $G=\GL(N,F)$. Here $\pi$ is an
irreducible supercuspidal representation of $G$.

Let $(J,\lambda)$ be a maximal simple type contained in $\pi$, as
in \cite{BK3}. We have $e=1$ and hence $\oA_\fs=\Mat(N,\integers_F)$. It
follows that $J^\fs=\GL(N,\mathfrak{o}_F)=L_0$.
By Proposition~\ref{simplet}, the representation
$\alpha=\Ind_J^{L_0}(\lambda)$ is irreducible.
The pair $(L_0,\alpha)$ is an $\fs$-type. The restriction to $L_0$ of a
smooth irreducible representation $\pi'$ of $G$ contains $\alpha$ if
and only if $\pi'$ is isomorphic to $\pi\otimes\chi\circ\det$,
where $\chi$ is an unramified quasicharacter of $F^{\times}$.
Moreover, $\pi$ contains $\alpha$ with multiplicity 1. In fact, the
representation $\alpha$ is the \emph{unique} smooth irreducible
representation $\tau$ of $L_0$ such that $(L_0,\tau)$ is an $\fs$-type,
see \cite{Pa}.

The little complex $C_*(\fs)$ determined by $\alpha$ is
\[
0 \longleftarrow C_0(\fs) \longleftarrow 0
\]
where $C_0 (\fs)$ is the free abelian group on the invariant
$0$-cycle
\[(\tau, \cR(\tau), \cR^2(\tau),\ldots,\cR^{n-1}(\tau))
\]

The total homology of the little complex is given by $h_0(\fs) =
\mathbb{Z}$. Therefore, by Lemma~\ref{Homology Lemma}, we have
\[H_{\ev}(\fs) = \mathbb{Z} = H_{\odd}(\fs).\]

\begin{thm}\label{supercuspidal}  Let $\pi$ be an irreducible unitary supercuspidal
representation of $\GL(N)$. Let $\fs = [G,\pi]_G$. Then we have
\[ H_{\ev}(\fs) \cong K_0(\fA(\fs)),\quad H_{\odd}(\fs) \cong K_1(\fA(\fs)).\]
\end{thm}
\begin{proof}The $C^*$-ideal $\fA(\fs)$ is given by
\[
\fA(\fs) = C(S^1, \mathfrak{K})
\] where $\mathfrak{K}$ is the $C^*$-algebra of compact operators and

\[S^1 =
\{ \pi \otimes \chi \circ \; {\rm det} \; : \chi \in
(F^{\times})^{\wedge} \}.\] The noncommutative $C^*$-algebra
$\fA(\fs)$ is strongly Morita equivalent to the commutative
$C^*$-algebra $C(S^1)$. For this $C^*$-algebra we have \[
K_j(C(S^1)) \cong K^j(S^1) = \mathbb{Z}\] where $j = 0,1$.

\end{proof}

\section{Generic Bernstein components attached to a maximal Levi subgroup}
We assume in this section that $\fs=[M,\sigma]_G$ with $M\cong
\GL(N_1)\times\GL(N_2)$ a $2$-blocks
Levi subgroup of $G$ such that $W_{\ft}=\{1\}$. Note that the last
conditions is always satisfied if $N_1\ne N_2$.

Let $(J_M,\lambda_M)$ be an $\ft$-type and let
$(J,\tau)$ be the $G$-cover of $(J_M,\tau_M)$ considered in
Theorem~\ref{ssinclusion}. We have shown there that
$\alpha:=\Ind_J^{J^\fs}(\tau)$ is irreducible.
It then follows from Propositions~\ref{IK} and \ref{Intbound} that
$\beta=\Ind_J^{L_0}(\tau)$ is irreducible.

Let $C_0 (\tau)$, $C_1(\tau)$ denote respectively the free abelian
group on one generator $(\beta, \cR(\beta), \ldots,\cR^{N-1}(\beta)$, and
on $(\alpha, \cR(\alpha),
\ldots, \cR^{N-1} \alpha)$. The little complex is
$$
0 \longleftarrow C_0(\fs) \stackrel{\partial}{\longleftarrow}
C_1(\fs) \longleftarrow 0
$$

\noindent The map $\partial$ is 0 by vertex compatibility of
$(\alpha, \cR(\alpha), \ldots, \cR^{N-1}(\alpha))$.  Then
$h_0(\fs) = \mathbb{Z}, h_1(\fs) = \mathbb{Z}$ and so
$H_{\ev}(\fs) = \mathbb{Z}^2 = H_{\odd}(\fs)$.

The subset of the tempered dual of $\GL(N)$ which contains the
$\fs$-type $(J, \tau)$ has the structure of a compact 2-torus. But
$K^0(\mathbb{T}^2) = \mathbb{Z}^2 = K^1 (\mathbb{T}^2)$ as
required.

\begin{thm} \label{special}
The $\fs$-type $(J,\tau)$ generates a little complex $C(\fs)$.
For this complex we have
\[
H_{\ev}(\fs) \cong K_0(\fA(\fs)) = \mathbb{Z}^2,\quad
H_{\odd}(\fs) \cong K_1(\fA(\fs)) = \mathbb{Z}^2\]
\end{thm}

Note that the above Theorem applies to the intermediate principal
series of $\GL(3)$. In the next section, we will consider the
principal series of $\GL(3)$.

\section{Principal series in $\GL(3)$}

Here $s_0$, $s_1$, $s_2$ are the standard involutions
\[ s_1 = \left( \begin{array}{ccc} 0 & 1 & 0\\ 1 & 0 & 0\\ 0 & 0
& 1
\end{array} \right) \qquad
s_2 = \left( \begin{array}{ccc} 1 & 0 & 0\\ 0 & 0 & 1\\ 0 & 1 & 0
\end{array} \right) \qquad
s_0 = \left( \begin{array}{ccc} 0 & 0 & \varpi^{-1}\\ 0 & 1 & 0\\
\varpi & 0 & 0
\end{array} \right)
\]
where
\[ \Pi =\Pi_3= \left( \begin{array}{ccc} 0 & 1 & 0\\ 0 & 0 & 1\\ \varpi
& 0 & 0
\end{array} \right).
\]

Note that $\val(\det(\Pi)) = 1$. Restricted to the affine line
$\mathbb{R}$ in the enlarged building $\beta^1\GL(3) = \beta
\SL(3) \times \mathbb{R}$, $\Pi$ sends $t$ to $t+1$. We also have
$\Pi^3 = \varpi 1 \in \GL(3)$.








We have the double coset identities

\begin{equation} \label{IJkI}
0 \leq k \leq 2 \Longrightarrow I \backslash J_k / I = \{1, s_k \}
\end{equation}
\begin{equation} \label{JrLsJr}
r,s,t \; {\rm distinct} \; \Longrightarrow J_r \backslash L_s /
J_r = \{1, s_t\}.
\end{equation}

Let $\fs =[T,\sigma]_G$, where $T$ is the diagonal split torus in
$\GL(3)$:
\[T=\left(\begin{matrix}
F^\times&0&0\cr
0&F^\times&0\cr
0&0&F^\times\cr\end{matrix}\right),\]
and $\sigma$ is an irreducible smooth character of $T$.

\smallskip
\noindent {\bf 9.1. Construction of an $\fs$-type, following
Roche}

For $u\in F$, we set
\[x_{1,2}(u)=\left(\begin{matrix}
1&u&0\cr
0&1&0\cr
0&0&1\cr\end{matrix}\right),\quad
x_{1,3}(u)=\left(\begin{matrix}
1&0&u\cr
0&1&0\cr
0&0&1\cr\end{matrix}\right),\quad
x_{2,3}(u)=\left(\begin{matrix}
1&0&0\cr
0&1&u\cr
0&0&1\cr\end{matrix}\right),\]
\[x_{2,1}(u)=\left(\begin{matrix}
1&0&0\cr
u&1&0\cr
0&0&1\cr\end{matrix}\right),\quad
x_{3,1}(u)=\left(\begin{matrix}
1&0&0\cr
0&1&0\cr
u&0&1\cr\end{matrix}\right),\quad
x_{3,2}(u)=\left(\begin{matrix}
1&0&0\cr
0&1&0\cr
0&u&1\cr\end{matrix}\right),\]
and, for any $k\in\ZZ$,
\[U_{i,j,k}=x_{i,j}(\fp_F^k).\]

\smallskip

Let $\Phi=\left\{\alpha_{i,j}\,:\,1\le i,j\le 2\right\}$ be the set of
roots of $G$ with respect to $T$.
For each root $\alpha_{i,j}$, let $\alpha_{i,j}^\vee$ denotes the
corresponding coroot.
We have
\[\alpha_{1,2}^\vee(t)=
\left(\begin{matrix}
t^{-1}&0&0\cr
0&t&0\cr
0&0&1\cr\end{matrix}\right),\quad
\alpha_{2,1}^\vee(t)=\left(\begin{matrix}
t&0&0\cr
0&t^{-1}&0\cr
0&0&1\cr\end{matrix}\right),\]
\[
\alpha_{1,3}^\vee(t)=\left(\begin{matrix}
t^{-1}&0&0\cr
0&1&0\cr
0&0&t\cr\end{matrix}\right),\quad
\alpha_{3,1}^\vee(t)=\left(\begin{matrix}
t&0&0\cr
0&1&0\cr
0&0&t^{-1}\cr\end{matrix}\right),\]
\[
\alpha_{2,3}^\vee(t)=\left(\begin{matrix} 1&0&0\cr 0&t^{-1}&0\cr
0&0&t\cr\end{matrix}\right),\quad
\alpha_{3,2}^\vee(t)=\left(\begin{matrix} 1&0&0\cr 0&t&0\cr
0&0&t^{-1}\cr\end{matrix}\right).\] Define $\sigma\colon
T\to\mathbb{T}$ by
\[\sigma\left(\begin{matrix}a&0&0\cr
0&b&0\cr
0&0&c\cr\end{matrix}\right)=\sigma_1(a)\sigma_2(b)\sigma_3(c),\]
where $\sigma_i\colon F^\times\to\mathbb{T}$ is a character of
$F^\times$, for $i=1,2,3$.

Hence
$\sigma\circ\alpha_{i,j}^\vee\colon\integers_F^\times\to\mathbb{T}$
is the smooth character of $\integers_F^\times$ defined by
\[\sigma\circ
\alpha_{i,j}^\vee(t)=\sigma_j(t)\sigma_i(t^{-1})=(\sigma_j\sigma_i^{-1})(t).\]

Now if $\chi\colon\integers_F^\times\to\mathbb{T}$ is a smooth
character, let $c(\chi)$ be the conductor of $\chi$: the least
integer $n\ge 1$ such that $1+\fp_F^n\subset\ker(\chi)$. We will
write $c_{i,j}$ for $c(\sigma\circ \alpha_{i,j}^\vee)$. We get
\[c_{i,j}=c(\sigma_j/\sigma_i)=c_{j,i}.\]
We can define a function $f=f_\sigma\colon\Phi\to\ZZ$ (here $\Phi$ is
the set of roots) as follows:
\[f_\sigma(\alpha_{i,j})=\begin{cases}[c_{i,j}/2]&\text{ if
$\alpha_{i,j}\in\Phi^+$,}\cr
[(c_{i,j}+1)/2]&\text{ if
$\alpha_{i,j}\in\Phi^-$.}\end{cases}\]
Here $[x]$ denotes the largest integer $\le x$.

Let
\[U_\sigma=\langle U_{i,j,f(\alpha_{i,j})}\,:\,\alpha_{i,j}\in\Phi\rangle,\]
and
\[J=\langle {}^\circ T,U_\sigma\rangle={}^\circ TU_\sigma=U_\sigma {}^\circ T,\]
where ${}^\circ T$ is the compact part of $T$,
\[{}^\circ T=\left(\begin{matrix}
\integers_F^\times&0&0\cr
0&\integers_F^\times&0\cr
0&0&\integers_F^\times\cr\end{matrix}\right).\]
It follows that
\[J=\left(\begin{matrix}
\integers_F^\times&\fp_F^{[c_{1,2}/2]}&\fp_F^{[c_{1,3}/2]}\cr
\fp_F^{[(c_{1,2}+1)/2]}&\integers_F^\times&\fp_F^{[c_{2,3}/2]}\cr
\fp_F^{[(c_{1,3}+1)/2]}&\fp_F^{[(c_{2,3}+1)/2]}&\integers_F^\times\cr
\end{matrix}\right).\]
The group $J$ will give the open compact group we are looking for.

\smallskip
Next, we need to figure out what is the correct character of $J$. In order
to do that, we set
\[T_{\sigma}=
\prod_{\alpha_{i,j}\in\Phi}\alpha_{i,j}^\vee(1+\fp_F^{f(\alpha_{i,j})+f(-\alpha_{i,j})})\,
\subset\,{}^\circ T.\]
Setting
\[U_\sigma^+=U_\sigma\cap\left(\begin{matrix}1&F&F\cr
0&1&F\cr
0&0&1\cr
\end{matrix}\right)\;\;\text{ and }\;\;
U_\sigma^-=U_\sigma\cap\left(\begin{matrix}
1&0&0\cr
F&1&0\cr
F&F&1\cr
\end{matrix}\right),\]
we obtain
\[U_\sigma=U_\sigma^-\cdot T_\sigma\cdot U_\sigma^+
\;\;\text{ and }\;\;
J=U_\sigma^-\cdot {}^\circ T\cdot U_\sigma^+.\]
It follows that
\[J/U_\sigma\cong {}^\circ T/T_{\sigma}.\]
By construction, $T_\sigma\subset\ker(\sigma_{|{}^\circ T})$.
Hence $\sigma_{|{}^\circ T}$ defines a character of ${}^\circ
T/T_{\sigma}$, and so can be lifted to a character $\tau$ of $J$.
Then $(J,\tau)$ an $\fs$-type by \cite[Theorem~7.7]{Roc}.

\smallskip
\noindent {\bf 9.2 Intertwining}

We first recall that the following results (\cite[Theorem~4.15]{Roc})
\begin{equation} \label{Itype}
I_G(\tau)=J\, \widetilde W(\sigma)\, J,
\end{equation}
where
\[\widetilde W(\sigma)=\left\{v\in\widetilde W\,:\,{}^v\sigma=\sigma\right\}.\]
More generally, it follows by the same proof as those of
\cite[Theorem~4.15]{Roc}, using \cite[Prop.~9.3]{AR} instead of
\cite[Prop.~4.11]{Roc}, that, for each $w\in W$,
\begin{equation} \label{Etypes}
I_G(\tau,{}^\tau)=J\, \widetilde W(\sigma,{}^w\sigma)\, {}^wJ,
\end{equation}
where
\[\widetilde W(\sigma,{}^w\sigma)=\left\{v\in\widetilde W\,:\,
{}^v\sigma={}^w\sigma \right\}.\]
Let
\[\Phi(\sigma)=\{\alpha_{i,j}\in\Phi\;:\;
(\sigma_i)_|{\integers_F^\times}=(\sigma_j)_|{\integers_F^\times}\}\,\subset
\,\Phi.\]
The group $W_0(\sigma)$ is equal to the group $W_{\fs_T}$, where
$\fs_T=[T,\lambda]_T$.
We observe that
\begin{equation} \label{compactI}
I_{L_0}(\tau)=J\,W_0(\sigma)\,J.
\end{equation}

\smallskip
\noindent {\bf 9.2.1 The case $\Phi(\sigma)=\Phi$}. Let $\fs
=[T,\sigma]_G$, where $\sigma = \psi \circ \det$ with $\psi$ a
smooth character of $F^{\times}$. In this case $c_{i,j}=1$ for any
$i$, $j$. It follows that $J=I$.

The pair $(I,\tau)$ is an $\fs$-type. We will construct cycles
from this type. It follows from~(\ref{compactI}) that, as
$\mathbb{C}$-algebras,
\[
\End_{L_0} \; (\Ind_I^{L_0}\tau)  \cong  {\cal H} (\GL (3,k_F)
// B).\] We also have, as $\mathbb{C}$-algebras, \[{\cal H} (\GL (3,k_F) // B)
\cong\mathbb{C} [W_0] \cong \mathbb{C} \oplus \mathbb{C} \oplus
M_2 (\mathbb{C})
\]
so that
$$
\Ind_I^{L_0}\tau = \lambda_{L_0} \oplus \mu_{L_0} \oplus \nu_{L_0}
\oplus \nu_{L_0}
$$
where $\lambda_{L_0}, \mu_{L_0},\nu_{L_0}$ are distinct.

We also have
\[
\sigma | {J_0} \hookrightarrow \Ind_I^{J_0}\tau
\]
by Frobenius reciprocity.  The triple $(\sigma | J_0,
\cR(\sigma|J_0),\cR^2(\sigma|J_0))$ is an invariant $1$-cycle, and
is not the boundary of $1_I$.

We now form the little complex:
\begin{itemize}
\item $C_0 (\fs)$ is the free abelian group on the three invariant
$0$-cycles
\[\lambda_L: = (\lambda_{L_0},\cR(\lambda_{L_0}),\cR^2(\lambda_{L_0}))\]
\[ \mu_L: = (\mu_{L_0}, \cR(\mu_{L_0}), \cR^2(\mu_{L_0}))\] \[\nu_L: =
 (\nu_{L_0}, \cR(\nu_{L_0}), \cR^2(\nu_{L_0})) \] \item $C_1(\fs)$ is the free abelian group
on the invariant $1$-cycle
\[\lambda_J: = (\sigma|J_0, \cR(\sigma|J_0), \cR^2(\sigma|J_0))\]
\end{itemize}
In the little complex\[ 0 \longleftarrow C_0(\fs)
\stackrel{0}{\longleftarrow} C_1(\fs)
\longleftarrow 0\]
we have
\[h_0(\fs) = \mathbb{Z}^3,\; h_1(\fs) = \mathbb{Z}.\]

The total homology of the little complex is $\mathbb{Z}^4$.
 As
generating cycles we may take \[ \lambda_L, \mu_L, \nu_L,
\lambda_J.\] and so, by Lemma 1, the even (resp. odd) chamber
homology groups are
\[H_{\ev}(\fs)= \mathbb{Z}^4, \;\; H_{\odd}(\fs) = \mathbb{Z}^4.\]

Each irreducible representation $\rho$ of a compact open subgroup
$J$ creates an \emph{idempotent} in $\cA$ as follows.  Let $d$
denote the dimension of $\rho$, let $\chi$ denote the character of
$\rho$.  Form the function $d \cdot \chi : J \longrightarrow
\mathbb{C}$ and \emph{extend by} $0$ to $G$. This function on $G$
is a non-zero idempotent in $\cA$, with the convolution product.
We will denote this idempotent by $e(\rho)$: \[ e(\rho) * e(\rho)
= e(\rho).\]

The inclusion
\[
H_{\ev}(\fs) \hookrightarrow K_0(\cA)
\]
is given explicitly as follows:
\[
\lambda_{L} \mapsto e(\lambda_{L_1}),\mu_L \mapsto
e(\mu_{L_1}),\nu_L \mapsto e(\nu_{L_1}), \lambda_J \mapsto
e(\lambda_{J_1}).\]

It follows from \cite{P} that the $C^*$-ideal $\fA(\fs)$ is given
as follows: \[ \fA(\fs) \cong C ( \; {\rm Sym}^3 \;
\mathbb{T},\mathfrak{K}) \oplus C (\mathbb{T}^2,\mathfrak{K})
\oplus C (\mathbb{T},\mathfrak{K}).\] The symmetric cube ${\rm
Sym}^3 \; \mathbb{T}$ is homotopy equivalent to $\mathbb{T}$ via
the product map \[{\rm Sym}^3 \; \mathbb{T} \sim
\mathbb{T},\;\;\;(z_1,z_2,z_3) \mapsto z_1z_2z_3.\] Hence $K_0
(\fA(\fs)) = \mathbb{Z}^4 = K_1(\fA(\fs))$ as required.

Note that
\begin{itemize}
\item ${\rm Sym}^3 \mathbb{T}$ is  in the minimal unitary
principal series of $\GL(3)$ \item $\mathbb{T}^2$ is in the
intermediate unitary principal series of $\GL(3)$ \item
$\mathbb{T}$ is in the discrete series of $\GL(3)$; if $\tau = 1$
then $\mathbb{T}$ comprises the unramified unitary twists of the
Steinberg representation of $\GL(3)$
\end{itemize}

These are precisely the tempered representations of $\GL(3)$ which
contain the type $(I,\tau)$.

\begin{thm} Let $\fs = [T,\sigma]_G$ where $\sigma = \psi \circ \det$ and $\psi$
is a smooth (unitary) character of $F^{\times}$.  Then we have \[
H_{\ev}(\fs) \cong K_0(\fA(\fs)) = \mathbb{Z}^4, \quad
H_{\odd}(\fs) \cong K_1(\fA(\fs)) = \mathbb{Z}^4.\]
\end{thm}


\smallskip
\noindent
{\bf 9.2.2 The case $\emptyset\ne \Phi(\sigma)\ne \Phi$}

Assume that $(\sigma_1)_{|\integers_F^\times}=(\sigma_2)_{|\integers_F^\times}
\ne(\sigma_3)_{|\integers_F^\times}$.
We have
\[J=\left(\begin{matrix}
\integers_F^\times&\integers_F&\fp_F^{[\ell/2]}\cr
\fp_F&\integers_F^\times&\fp_F^{[\ell/2]}\cr
\fp_F^{[(\ell+1)/2]}&\fp_F^{[(\ell+1)/2]}&\integers_F^\times\cr
\end{matrix}\right),\]
where $\ell=c_{1,3}=c_{2,3}$, and
$$ \tau \left( \begin{array}{ccc} a & \ast & \ast \\
\ast & b & \ast \\ \ast & \ast & c \end{array} \right) =
\sigma_1(a)\sigma_1(b)\sigma_3(c).$$

It is clear that $s_1 \in I_{L_0}(\tau)$. The Weyl group
$W_{\fs_T} = \mathbb{Z}/2\mathbb{Z}$ and so we have
$I_{L_o}(\tau) = J \cup Js_1J$. The complete list is as follows:
\[
I_I(\tau) = J\]
\[
I_{J_1}(\tau) = J<s_1>J,\quad I_{J_2}(\tau) = J,\quad
I_{J_0}(\tau) = J\]
\[
I_{L_1}(\tau) = J<s'>J,\quad I_{L_2}(\tau) = J<s_1>J,\quad
I_{L_0}(\tau) = J<s_1>J\]
where
$$ s' = \left(
\begin{array}{lll} 0 & \varpi^{-1} & 0\\ \varpi & 0 & 0\\ 0 & 0 & 1
\end{array}\right)$$

\begin{lem}  Let $\tau_{1} = \Ind_{J}^{I}
(\tau)$. Then $\tau_1$ is irreducible.
\end{lem}
\begin{proof} This follows from proposition~\ref{IK}, since $I_I(\tau) = J$.
It follows that $(I, \tau_1)$ is an $\fs$-type.\end{proof}

\begin{lem} We have
\[
{\Ind}^{J_1}_{I} \tau_1  =  \xi_1 \oplus \eta_1,\quad
{\Ind}^{L_0}_{I} \tau_1  =  \gamma_0 \oplus \delta_0.
\]
\end{lem}
\begin{proof}  We have $I_{J_1}(\tau) = J \cup s_1J$. Hence
$$
\End_{J_{1}} (\Ind^{J_{1}}_{I} \tau_1) =  \cI_1 (\tau) \oplus
\cI_{s_{1}} (\tau) = \mathbb{C} \oplus \mathbb{C}.
$$
This implies that $\Ind^{J_{1}}_{I} \tau_1$ has two distinct
irreducible constituents $\xi_1$, $\eta_1$.  Now, we replace $J_1$
by $L_0$, and infer that $\Ind^{L_0}_{I} \tau_1$ has two distinct
irreducible constituents $\gamma_0, \delta_0$.
\end{proof}

It follows that
\[
\Ind_I^{J_2} \,\cR(\tau_1)  = \cR(\xi_1) \oplus \cR(\eta_1),\]
\[
\Ind_I^{J_0} \,\cR^2(\tau_1) = \cR^2(\xi_1) \oplus
\cR^2(\eta_1).\]

This creates two invariant $1$-chains

\[\xi: = (\xi_1,\cR(\xi_1),\cR^2(\xi_1)),
\;\eta: = (\eta_1, \cR(\eta_1), \cR^2(\eta_1)) .\]

It follows from~(\ref{RI}) that
\[\Ind_I^{L_0}\tau_1 \cong \Ind_I^{L_0}\cR(\tau_1)\] \[ \zeta_1: =
\Ind_I^{J_1}\cR(\tau_1) \cong \Ind_I^{J_1}\cR^2(\tau_1)\] By (17)
we have \[0 = \langle \Ind_J^{J_1}\tau,
\Ind_J^{J_1}\cR(\tau)\rangle.\]

Let $C_0 (\fs)$ be the free abelian group generated by the two
invariant $0$-cycles
\[(\gamma_0, \cR(\gamma_0), \cR^2(\gamma_0)), \;
(\delta_0, \cR(\delta_0), \cR^2(\delta_0)) .\] Let $C_1 (\fs)$ be
the free abelian group generated by the two invariant $1$-cycles
$\xi$ and $\zeta$.

The little complex is then $$ 0 \longleftarrow C_0(\fs)
\stackrel{0}{ \longleftarrow} C_1(\fs)
\longleftarrow 0.$$

We have $h_0(\fs) = \mathbb{Z}^2, h_1(\fs) = \mathbb{Z}^2$ and the
total homology is $\mathbb{Z}^4$ and so $H_{\ev}(\fs) =
\mathbb{Z}^4 = H_{\odd}(\fs)$.

The definition of $\zeta: = (\zeta_0,\cR(\zeta_0),\cR^2(\zeta_0))$
shows that
\[
\partial(\tau_1 + \cR(\tau_1) + \cR^2(\tau_1)) = \xi + \eta + 2\zeta
\]
so that $\eta$ and $-(\xi + 2\zeta)$ are homologous. Therefore the
invariant $1$-cycle $\eta$ does not contribute a new homology
class in $H_1(G;\beta^1G)$.

The $C^*$-ideal $\fA(\fs)$ is as follows: \[
C(\mathbb{T}^2,\mathfrak{K})\oplus C(\Sym^2\mathbb{T} \times
\mathbb{T},\mathfrak{K}).\]

To identify these ideals, we proceed as follows. First, let
$\Psi(F^{\times})$ denote the group of unramified unitary
characters of $F^{\times}$. The first summand is determined by the
compact orbit
\[\mathcal{O}(\St(\sigma_1,2) \otimes \sigma_3) =
\{\chi_1\St(\sigma_1,2) \otimes \chi_2\sigma_3: \chi_j \in
\Psi(F^{\times})\}\] where $\St(\sigma_1,2)$ is a generalized
Steinberg representation; the second is determined by the compact
orbit
\[\mathcal{O}(\sigma_1 \otimes \sigma_1 \otimes \sigma_3) =
 \{\chi_1 \sigma_1 \otimes \chi_2\sigma_1 \otimes
\chi_3\sigma_3: \chi_j \in \Psi(F^{\times})\}.\]  The compact
space $\Sym^2 \mathbb{T} \times \mathbb{T}$  is homotopy
equivalent to the $2$-torus $\mathbb{T}^2$.

The space $\Sym^2\mathbb{T} \times \mathbb{T}$ is in the minimal
unitary principal series of $\GL(3)$ and the space $\mathbb{T}^2$
is in the intermediate unitary principal series of $\GL(3)$. The
union of these two compact spaces is precisely the set of tempered
representations of $\GL(3)$ which contain the $\fs$-type
$(J,\tau)$.

The $K$-groups are now immediate: \[ K_j(\fA(\fs)) = \Zset^4\]with
$j = 0,1$.

\begin{thm}  Let $\fs = [T,\sigma]_G$.  We have
\[
H_{\ev}(\fs) \cong K_0(\fA(\fs)) = \mathbb{Z}^4, \quad
H_{\odd}(\fs) \cong K_1(\fA(\fs)) = \mathbb{Z}^4.\]
\end{thm}

\smallskip
\noindent {\bf 9.2.3 The case $\Phi(\sigma)=\emptyset$}. The
generic torus.  The Bernstein component is $[T,\sigma_1 \otimes
\sigma_2 \otimes \sigma_3]$.  The Weyl group $W(T) = W_0 =S_3$,
and the associated parahoric subgroup is the Iwahori subgroup $I$.

The restrictions of $\sigma_1$, $\sigma_2$ and $\sigma_3$ to
$\integers_F^\times$ are all distinct. We have
$\Phi(\sigma)=\emptyset$.
We have $\widetilde W(\sigma)=D$, where $D$ is the subgroup of $T$ whose
eigenvalues are powers of $\varpi$.  The subgroup $D$ is free abelian of rank 3.
The only compact element in $D$ is $1_G$.
The only double-$J$-coset representative in
$L_0$ which $G$-intertwines $\tau$ is $1_G$.  This proves the
following:

\begin{lem}  If $r = 0,1,2$ then $\Ind^{J_r}_{J} (\tau)$ is irreducible,
$\Ind^{L_r}_{J}(\tau)$ is irreducible.
\end{lem}

Let $\alpha = \Ind^I_J(\tau)$. Then $\alpha$ is irreducible.
Therefore $(I,\alpha)$ is an $\fs$-type.

\begin{lem} If $w \in W_0$ then ${\rm
Ind}^{L_0}_{I} \alpha = \Ind^{L_0}_{I} ({}^w\alpha)$.
\end{lem}
\begin{proof}  We have
$\Ind_I^{L_0}(\alpha)=\Ind_J^{L_0}(\tau)$ and
$\Ind_I^{L_0}({}^w\alpha)=\Ind_J^{L_0}({}^w\tau)$.
By Proposition~\ref{IKequiv}, it is sufficient to prove that
$I_G(\tau,{}^w\tau)\ne\{0\}$. But $I_G(\tau,{}^w\tau)=J\,\widetilde
W(\sigma,{}^w\sigma)\,J$.
\end{proof}

\begin{lem} \label{invariant}
If $w \in W_0$ then
\[\Ind^{J_r}_{I} (\alpha) \cong \Ind^{J_r}_{I} (^w \alpha)
\Longleftrightarrow w \in <s_r>\]
with $0 \leq r \leq 2$.
\end{lem}
\begin{proof} By Proposition~\ref{IKequiv},
\[\Ind^{J_r}_{J} (\tau) \cong \Ind^{J_r}_{J} (^w \tau)
\Longleftrightarrow
I_{J_r}(\tau,{}^w\tau)\ne\{0\}.\]
From~(\ref{Etypes}), we have
\[I_{J_r}(\tau,{}^w\tau)=J_r\,\cap\,\widetilde W(\sigma,{}^w\sigma)=
J_r\,\cap\,\widetilde W(\sigma)\cdot w=
J_r\,\cap\,D\cdot w.\]
The result follows from the fact that $J_r=I\,<1,s_r>\,I$.
\end{proof}

Inducing the orbit $W_0 \cdot \alpha$ from $J$ to $J_1$ gives 3
distinct elements $\rho_1, \phi_1, \psi_1$, by Lemma 9.  Inducing
from $J$ to $L_0$ gives $\gamma_0$.

 Set $C_2 (\fs)$ = free abelian group on the invariant $2$-cycle
 \[\epsilon: = \sum_{w \in W_0}
 sgn(w)(^w \alpha).\]

Set $C_1 (\fs)$ = free abelian group on the three invariant
$1$-cycles
\[\rho:=(\rho_1, \cR(\rho_1), \cR^2(\rho_1)),\] \[\phi:= (\phi_1, \cR(\phi_1),
\cR^2(\phi_1)),\] \[\psi: = (\psi_1, \cR(\psi_1),
\cR^2(\psi_1)).\]

Set $C_0 (\fs)$ = free vector abelian group on the invariant
$0$-cycle
\[\gamma: = (\gamma_0, \cR(\gamma_0), \cR^2(\gamma_0)).\]

Note that
\[
\partial (\sum_{w \in Alt(3)} \,^w \alpha) = \rho +  \phi + \psi
\]
where $Alt(3)$ is the alternating subgroup of $W_0$.  Since
$^{s_1s_2} \alpha = \cR(\alpha)$, we may also write this as
\[
\partial(\alpha + \cR(\alpha) + \cR^2(\alpha)) = \rho +  \phi +
\psi.
\]
It follows that $\psi$ is homologous to $- (\rho + \phi)$ in the
top row of the double complex $C_{**}$. This implies that the
\emph{image} of $C(\fs)$ in $C_{**}$ determines $4$ homology
classes. As representing cycles we may take the $2$-cycle
$\epsilon$, the two $1$-cycles $\rho, \phi$, and the $0$-cycle
$\gamma$. Therefore
\[\RH_{\ev}(G;\beta^1G)^{\fs} = \RH_{\odd}(G;\beta^1G)^{\fs} =
\mathbb{Z}^4.\]

\begin{thm}The subspace of the tempered dual of $\GL(3)$ which contains the
$\fs$-type $(I, \alpha)$ has the structure of a compact 3-torus.
This is a generic torus in the minimal unitary principal series of
$\GL(3)$. We have
\[\RH_{\ev}(G;\beta^1G)^{\fs} \cong K_0(\fA(\fs)) = \Zset^4,
\quad \RH_{\odd}(G;\beta^1G)^{\fs} \cong K_1(\fA(\fs)) =
\Zset^4.\]
\end{thm}
\begin{proof}  Let $\Psi(F^{\times})$ denote the group of unramified characters of $F^{\times}$.
If $\chi \in \Psi(F^{\times})$ then $\chi(x) = z^{val(x)}$ with
$z$ a complex number of modulus $1$, so that
\[ \Psi(F^{\times}) \cong \mathbb{T}.\] Writing
\[ \mathbb{T}^3 = \{\Ind_T^G(\chi_1\sigma_1 \otimes \chi_2\sigma_2
\otimes \chi_3 \sigma_3): \chi_j \in \Psi(F^{\times})\}\] we have
\[\fA(\fs) \cong C(\mathbb{T}^3, \mathfrak{K})\] which is strongly
Morita equivalent to $C(\mathbb{T}^3)$. The $K$-theory of the
$3$-torus is given by
\[
K^j(\mathbb{T}^3) = \mathbb{Z}^4\] where $j = 0,1$.
\end{proof}

\appendix
\section{Chamber homology and K-theory}

Let $G = \GL(N)$ and let $\cA$ denote the reduced $C^*$-algebra of
$G$. Let $\mathcal{H}(G)$ be the convolution algebra of uniformly
locally constant, compactly supported, complex-valued functions on
$G$, and let $\mathcal{C}(G)$ be the Harish-Chandra Schwartz
algebra of $G$. The following diagram serves as a framework for
this article:
$$
\begin{diagram}
\node{K_j^{\text{top}}(G)}  \arrow{s,l}{\text{ch}}
\arrow[2]{e,t}{\mu}   \node[2]{K_j(\cA)}
\arrow{s,r}{\text{ch}}  \\
  \node{\RH_j(G; \beta^1 G)\otimes_{\Zset} \Cset} \arrow{e}
 \node{ \HP_j(\mathcal{H}(G))} \arrow{e,t}{\imath_*}
 \node{\HP_j(\mathcal{C}(G))}\\
\end{diagram}
$$ with $j = 0,1$.  In this diagram,  $K_j^{\text{top}}(G)$ denotes the
topological $K$-theory of $G$, $K_j(\cA)$ denotes $K$-theory for
the $C^*$-algebra $\cA$. In addition, $\HP_j(\mathcal{H}(G))$
denotes periodic cyclic homology of the algebra $\mathcal{H}(G)$,
and $\HP_j(\mathcal{C}(G))$ denotes periodic cyclic homology of
the topological algebra $\mathcal{C}(G)$.  For periodic cyclic
homology, see \cite[2.4]{CST}.

The Baum-Connes assembly map $\mu$ is an isomorphism
\cite{BHP1,L}. The map
$$
\RH_*(G; \beta^1 G)\otimes_{\Zset} \Cset \longrightarrow
\HP_*(\mathcal{H}(G))
$$
is an isomorphism \cite{HN,S}. The map $\imath_*$ is an
isomorphism by \cite{BHP1,BP1}. The right hand Chern character is
constructed  in \cite{BP2} and is an isomorphism after tensoring
over $\mathbb{Z}$ with $\bbc$ \cite[Theorem 3]{BP2}. The left hand
Chern character is the unique map for which the diagram is
commutative.

\section{The Bernstein spectrum}

Let $G$ be the group of $F$-points of a connected reductive
algebraic group defined over $F$.  We consider pairs $(L,\sigma)$
where $L$ is a Levi subgroup of a parabolic subgroup of $G$, and
$\sigma$ is an irreducible
supercuspidal representation of $L$.  We say two such pairs $(L_1,\sigma_1)$,
$(L_2,\sigma_2)$ are \emph{inertially
equivalent} if there exist $g \in G$ and an unramified
quasicharacter $\chi$ of $L_2$ such that \[ L_2 = L_1^g \quad
\text{and} \quad \sigma_1^g \cong \sigma_2 \otimes \chi.\] Here,
$L_1^g: = g^{-1}L_1g$ and $\sigma_1^g(x) = \sigma_1(gxg^{-1})$ for
all $x \in L_1^g$.  We write $[L,\sigma]_G$ for the inertial
equivalence of the pair $(L,\sigma)$ and $\fB(G)$ for the set of
all inertial equivalence classes.   The set $\fB(G)$ is the
\emph{Bernstein spectrum} of $G$.  We will write $\fs \in \fB(G)$.

The Hecke algebra $\cH(G)$ is a unital $\cH(G)$-module via left
multiplication, and admits the canonical Bernstein decomposition
as a purely algebraic direct sum of two-sided ideals:
\[
\cH(G)= \bigoplus_{\fs \in \mathfrak{B}(G)}\cH(G)^{\fs}.\] This
determines the canonical Bernstein decomposition of the reduced
$C^*$-algebra as a $C^*$-direct-sum of two-sided $C^*$-ideals:
\[ \cA = \bigoplus_{\fs \in \mathfrak{B}(G)}\fA(\fs).\]

Now $C^*$-direct sums are respected by the $K$-theory of
$C^*$-algebras, and we have
\begin{eqnarray} K_j(\cA) = \bigoplus_{\fs \in \mathfrak{B}(G)}
K_j(\fA(\fs))\end{eqnarray} with $j = 0,1$. The abelian groups
$K_j(\fA(\fs))$ are finitely generated free abelian groups, see
\cite{P}.

We will define $\RH_{\ev/\odd}(G; \beta^1 G)^{\fs}$ as the
pre-image of $K_j(\fA(\fs))$ via the commutative diagram in
Appendix A:
\begin{eqnarray} \RH_{\ev}(G; \beta^1 G)^{\fs} \cong
K_0(\fA(\fs)), \quad \RH_{\odd}(G; \beta^1 G)^{\fs} \cong
K_1(\fA(\fs)).\end{eqnarray}

\section{The formula for the rank}

Let $\fs$ be a point in the Bernstein spectrum $\fB(G)$, so that
$\fs = [L,\sigma]_G$.   We have \[ L = \GL(m_1)^{e_1} \times
\cdots \times \GL(m_r)^{e_r}\] with $m_1e_1 + \cdots + m_re_r =
N$. The numbers $e_1, \ldots,e_r$ are called the \emph{exponents}
of $\fs$, as in \cite{BP1}. According to \cite[Lemma 3.2]{BP1}, we
then have
\begin{eqnarray} \textrm{rank}\, K_j(\fA(\fs))
= 2^{r-1} \beta(e_1) \cdots \beta(e_r)
\end{eqnarray} where \[ \beta(e) = \sum
2^{\kappa(\pi) - 1}.
\]
In this formula, $\pi$ is a partition of $e$, the sum is over all
partitions of $e$, and $\kappa(\pi)$ is the number of unequal
parts of $\pi$.  For example, if $\pi$ is the partition
$1+1+1+3+3+3+3+7+9$ of $31$ then $\kappa(\pi) = 4$.

The ranks of the finitely generated abelian groups
$\RH_{\ev/odd}(G; \beta^1 G)^{\fs}$ are given by
\begin{eqnarray}\textrm{rank} \, \RH_{\ev/\odd}(G; \beta^1 G)^{\fs} =  2^{r-1} \beta(e_1) \cdots
\beta(e_r).
\end{eqnarray}

\section{Invariants attached to $\fs$}

We write the supercuspidal representation $\sigma$ of the Levi
subgroup
$$M\,\cong\,\prod_{i=1}^q\prod_{j=1}^{c_i}\GL(N_{i,j},F)$$ as a vector
$\sigma=(\sigma_{1,1}, \ldots, \sigma_{1,{c_1}},
\sigma_{2,1},\ldots,\sigma_{2,{c_2}}, \ldots, \sigma_{q,1},
\ldots,\sigma_{q,{c_q}})$ where $\sigma_{i,j}$ is an irreducible
supercuspidal representation of $\GL(N_{i,j},F)$, and for each
$i\in\{1,\ldots,q\}$, the representations $\sigma_{i,j}$ ($1\le
j\le c_i$) admit the same endo-class. At the same time, for all
$1\le j\le c_i$ and $1\le j'\le c_{i'}$, the representations
$\sigma_{i,j}$ and $\sigma_{{i'},{j'}}$ have distinct endo-classes
if $i'\ne i$.  This implies that, for a given $i$, in the
construction of Bushnell-Kutzko, all the representations
$\sigma_{i,j}$ ($1\le j\le c_i$) may be assumed to correspond to the same field
extension $E_i$ of $F$. Let $e(E_i|F)$ denote the ramification index of
$E_i$ over $F$. Then the parahoric subgroup $J^{\fs}$ only
depends on the integers $N_{i,j}$, $c_i$ and $e(E_i|F)$ (see
Definition~\ref{Jsgeneral}).

For supercuspidal representations, the parahoric subgroup is
always the same one, say $\GL(N,\integers_F)$;  when $q=1$ (that
is, only one endo-class), the parahoric is given by the integers
$N_{1,1}$, ..., $N_{1,c_1}$, which are the sizes of the blocks of
$M$. In the general case, the parahoric subgroup depends on the
sizes of the blocks of $M$, of the block decomposition defined by
the endo-classes (that is, those corresponding to the Levi
subgroup $\bar M\cong\prod_{i=1}^q\GL(\barN_i)$, with
$\barN_i=\sum_{j=1}^{c_i}N_{i,j}$) and on the ramification
indices.

Anne-Marie Aubert, Institut de Math\'ematiques de Jussieu,
U.M.R. 7586 du C.N.R.S., 175 rue du Chevaleret 75013 Paris, France.\\
Email: aubert@math.jussieu.fr

Samir Hasan, Department of Pure Mathematics, Faculty of Sciences,
University of Damascus, Damascus, S.A.R., SYRIA.\\
Email: samir.hasan@gmail.com

Roger Plymen, School of Mathematics, Manchester University, M13
9PL, England.\\
Email: plymen@manchester.ac.uk

\begin{thebibliography}{99}
\bibitem{AR} Adler, J. and Roche, A.: An intertwining result for $p$-adic
groups, \emph{Canad. J. Math.} {\bf 52} (2000) 449--467.
\bibitem{BCH} Baum, P.,  Connes, A. and Higson, N.: Classifying space for proper
actions and $K$-theory of group $C^*$-algebras, \emph{Contemporary
Math.} {\bf 167} (1994) 241--291.
\bibitem{BHP1}  Baum, P., Higson, N. and Plymen, R.J.:
A proof of the Baum-Connes conjecture for $p$-adic $\GL(n)$,
\emph{C. R. Acad. Sci. Paris} {\bf 325} (1997) 171-176.
\bibitem{BHP2} Baum, P., Higson, N. and Plymen, R.J.:  Representations of
$p$-adic groups: a view from operator algebras, \emph{Proc. Symp.
Pure Math.} {\bf 68} (2001) 111--149.
\bibitem{Be} Bernstein, J. (r\'e\-di\-g\'e par P. Deligne): Le ``centre''
de Bernstein. \emph{Re\-pr\'e\-sen\-ta\-tions des grou\-pes
r\'e\-duc\-tifs sur un corps local.} Hermann, Paris (1984) 1--32.
\bibitem{BP1} Brodzki, J. and  Plymen, R.J.:  Complex structure on the smooth dual
of $\GL(n)$, \emph{Documenta Math.} {\bf 7} (2002) 91--112.
\bibitem{BP2} Brodzki, J. and Plymen, R.J.: Chern character for the Schwartz
algebra of $p$-adic $\GL(n)$, \emph{Bull. London Math. Soc.} {\bf
34} (2002) 219--228.
\bibitem{BK3} Bushnell, C.J. and Kutzko, P.C.: The admissible dual of $\GL(N)$
via compact open subgroups, \emph{Annals of Math. Study} {\bf 129}
(1993) Princeton University Press.
\bibitem{BK1} Bushnell, C.J. and  Kutzko, P.C.:
Smooth representations of reductive $p$-adic groups: structure
theory via types, \emph{Proc. London Math. Soc.} {\bf 77} (1998)
582--634.
\bibitem{BK2} Bushnell, C.J. and Kutzko, P.C.:
Semisimple types in $\GL_n$, \emph{Compos. Math.} {\bf 119} (1999)
53--97.
\bibitem{BK4} Bushnell, C.J. and  Kutzko, P.C.: Types in reductive $p$-adic
groups: the Hecke algebra of a cover, \emph{Proc. Amer. Math.
Soc.} {\bf 129} (2001) 601--607.
\bibitem{CST} Cuntz, J., Skandalis, G. and Tsygan, B.: \emph{Cyclic homology
in noncommutative geometry}, EMS 121, Springer-Verlag, Berlin
2004.
\bibitem{GM} Gelfand, S.I. and Manin, Yu. I.: \emph{Homological algebra}, Springer-Verlag, Berlin 1999.
\bibitem{HN} Higson, N. and Nistor, V.:  Cyclic homology of totally disconnected
groups acting on buildings, \emph{J. Functional Analysis} {\bf
141} (1996) 466--485.
\bibitem{L} Lafforgue, V.: $K$-th\'{e}orie bivariante pour les alg\`{e}bres de
Banach et conjecture de Baum-Connes, \emph{Invent. Math.} {\bf
149} (2002) 1--95.
\bibitem{Pa} Paskunas, V.: Unicity of types for supercuspidal representations of
$\GL_N$, \emph{Proc. London Math. Soc.} {\bf 91}  (2005) 623--654.
\bibitem{P} Plymen, R.J.: Reduced $C^*$-algebra of the $p$-adic
group $\GL(n)$, \emph{J. Functional Analysis} {\bf 72} (1987)
1--12.
\bibitem{Roc} Roche, A.: Types and Hecke algebras for principal
series representations of split reductive $p$-adic groups,
\emph{Ann. scient. \'Ec. Norm. Sup.} {\bf 31} (1998) 361--413.
\bibitem{Ron} Ronan, M.: {\em Lectures on buildings}, Academic press
(1989).
\bibitem{S} Schneider, P.: Equivariant homology for totally
disconnected groups, \emph{J.Algebra} {\bf 203} (1998) 50--68.
\bibitem{T} Tits, J.: Reductive groups over local fields, in \emph{
Automorphic forms, representations and $L$-functions}, \emph{Proc.
Symp. Pure Math.} {\bf 33} (1979), part 1, 29--69.
\end{thebibliography}
\end{document}